\documentclass{amsart}
\usepackage{amsfonts}
\usepackage{comment}
\usepackage{lineno}
\newtheorem{theorem}{Theorem}
\newtheorem{lemma}[theorem]{Lemma}
\numberwithin{theorem}{section} 
\newtheorem{corollary}[theorem]{Corollary}
\newtheorem{definition}[theorem]{Definition}

\newcounter{lemmacounter}
\usepackage{graphicx}
\usepackage[dvips]{color}
\usepackage{amssymb}
\usepackage{listings}
\usepackage{pstricks}
\psset{unit=3cm}
\usepackage{url}
\usepackage[section]{placeins}	
\usepackage{floatrow}
\usepackage{caption}
\usepackage{alltt}
\usepackage{comment}
\usepackage[greek,english]{babel}
\usepackage[utf8x]{inputenc}
\newrgbcolor{lightblue}{0.8 0.8 1}
\newrgbcolor{pink}{1 0.8 0.8}
\newrgbcolor{lightgreen}{0.8 1 0.8}
\newrgbcolor{lightyellow}{1 1 0.8}

\def\B{{\bf B}}   

\def\FigureAdditiveOne{%
\psset{unit=1cm}
\pspicture(-0.4,-0.4)(4.2,3.2)
\pspolygon[fillcolor=lightblue,fillstyle=solid](0,0)(0.5,2) (2,0)
\pspolygon[fillcolor=lightyellow,fillstyle=solid](0.5,2)(3.5,3)(4,0)(2,0)
\psline[linestyle=dashed](0.5,2)(4,0)
\psdot(0,0)
\put(-0.25,-0.35){$B$}
\psdot(0.5,2)
\put(0.15,1.9){$A$}
\psdot(2,0)
\put(1.9,-0.35){$C$}
\psdot(3.5,3)
\put(3.6,2.9){$E$}
\psdot(4,0)
\put(3.9,-0.35){$D$}
\endpspicture
}

\def\FigureAdditiveTwo{%
\psset{unit=1cm}
\pspicture(-0.4,-0.4)(4.2,3.2)
\pspolygon[fillcolor=lightblue,fillstyle=solid](0,0)(0.5,2)(2,2.5)(2,0)
\pspolygon[fillcolor=lightyellow,fillstyle=solid](2,2.5)(3.5,3)(4,0)(2,0)
\psline[linestyle=dashed](0.5,2)(2,0)
\psdot(0,0)
\put(-0.25,-0.35){$F$}
\psdot(0.5,2)
\put(0.12,1.9){$A$}
\psdot(2,0)
\put(1.9,-0.35){$E$}
\psdot(3.5,3)
\put(3.6,2.9){$C$}
\psdot(4,0)
\put(3.9,-0.35){$D$}
\psdot(2,2.5)
\put(1.9,2.65){$B$}
\endpspicture
}

\def\FigurePasteThree{%
\psset{unit=1cm}
\pspicture(-0.4,-0.35)(9.5,3.4)
\pspolygon[fillcolor=lightblue,fillstyle=solid](0,1)(1,3)(1,0)
\pspolygon[fillcolor=lightyellow,fillstyle=solid](1,3)(1,0)(3,2)
\psline[linestyle=dashed](0,1)(3,2)
\psdot(0,1)
\put(-0.4,0.9){$A$}
\psdot(1,3)
\put(0.8,3.12){$B$}
\psdot(3,2)
\put(3.1,1.95){$C$}
\psdot(1,0)
\put(0.8,-0.35){$D$}
\psdot(1,1.33)
\put(1.13,1.05){$M$}
\pspolygon[fillcolor=lightblue,fillstyle=solid](4.5,2)(6,2)(6,0)
\pspolygon[fillcolor=lightyellow,fillstyle=solid](9,0.6)(6,2)(6,0)
\psline[linestyle=dashed](9,0.6)(4.5,2)
\psdot(4.5,2)
\put(4.15,1.9){$a$}
\psdot(6,2)
\put(5.9,2.12){$b$}
\psdot(9,0.6)
\put(9.12,0.53){$c$}
\psdot(6,0)
\put(5.9,-0.33){$d$}
\psdot(6,1.537)
\put(6.1,1.2){$m$}
\endpspicture
}

\def\FigurePasteThreeVerify{%
\psset{unit=1cm}
\pspicture(-0.4,-0.35)(9.5,3.4)
\pspolygon[fillcolor=lightblue,fillstyle=solid](0,1)(1,3)(1,0)
\pspolygon[fillcolor=lightyellow,fillstyle=solid](1,3)(1,0)(3,2)
\psline[linestyle=dashed](0,1)(3,2)
\pspolygon[linestyle=solid,linecolor=red,fillstyle=none]
(0,0)(0,3)(3,3)(3,0)
\pspolygon[linestyle=solid,linecolor=red,fillstyle=none]
(4.5,0)(4.5,2)(9,2)(9,0)
\psdot(0,1)
\put(-0.4,0.9){$A$}
\psdot(1,3)
\put(0.8,3.12){$B$}
\psdot(3,2)
\put(3.1,1.95){$C$}
\psdot(1,0)
\put(0.8,-0.35){$D$}
\psdot(1,1.33)
\put(1.13,1.05){$M$}
\pspolygon[fillcolor=lightblue,fillstyle=solid](4.5,2)(6,2)(6,0)
\pspolygon[fillcolor=lightyellow,fillstyle=solid](9,0.6)(6,2)(6,0)
\psline[linestyle=dashed](9,0.6)(4.5,2)
\psdot(4.5,2)
\put(4.15,1.9){$a$}
\psdot(6,2)
\put(5.9,2.12){$b$}
\psdot(9,0.6)
\put(9.12,0.53){$c$}
\psdot(6,0)
\put(5.9,-0.33){$d$}
\psdot(6,1.537)
\put(6.1,1.2){$m$}
\endpspicture
}

\def\FigureFLMKRectangle{%
\psset{unit=1cm}
\pspicture(-1.3,-0.25)(2,3)
\pspolygon[fillcolor=lightblue,fillstyle=solid](0.8,0)(0.8,1.5)(2,1.5)(2,0)
\pspolygon[fillcolor=lightyellow,fillstyle=solid](0.8,0)(0,0)(0,1.5)(0.8,1.5)
\psdot(0,0)
\put(-0.4,-0.35){$K$}
\psdot(2,0)
\put(2.1,-0.35){$M$}
\psdot(0,1.5)
\put(-0.4,1.4){$F$}
\psdot(2,1.5)
\put(2.1,1.4){$L$}
\psdot(0.8,0)
\put(0.6,-0.35){$H$}
\psdot(0.8,1.5)
\put(0.6,1.65){$G$}
\endpspicture
}

\def\FigureFLMKBigRectangle{%
\psset{unit=1cm}
\pspicture(-1.3,-0.25)(2,3.2)
\pspolygon(0,0)(2,0)(2,3)(0,3)
\pspolygon[fillcolor=lightblue,fillstyle=solid](0.8,0)(0.8,1.5)(2,1.5)(2,0)
\pspolygon[fillcolor=lightyellow,fillstyle=solid](0.8,0)(0,0)(0,1.5)(0.8,1.5)
\psdot(0,0)
\put(-0.4,-0.35){$K$}
\psdot(2,0)
\put(2.1,-0.35){$M$}
\psdot(0,1.5)
\put(-0.4,1.4){$F$}
\psdot(2,1.5)
\put(2.1,1.4){$L$}
\psdot(0.8,0)
\put(0.6,-0.35){$H$}
\psdot(0.8,1.5)
\put(0.6,1.65){$G$}
\psdot(0,3)
\put(-0.4,2.9){$k$}
\psdot(0.8,3)
\put(0.6,3.15){$h$}
\psdot(2,3)
\put(2.1,2.9){$m$}
\endpspicture
}

\def\FigureArchimedesOne{%
\psset{unit=1cm}
\pspicture(-1,-1)(1,1)
\pscircle(0.00,0.00){1.000000}
\pswedge[linestyle = none, fillcolor = white, fillstyle = solid]
(0.000000,0.000000){1.000000}{0.000000}{15.000000}
\pswedge[linestyle = none, fillcolor = black, fillstyle = solid]
(0.000000,0.000000){1.000000}{15.000000}{30.000000}
\pswedge[linestyle = none, fillcolor = white, fillstyle = solid]
(0.000000,0.000000){1.000000}{30.000000}{45.000000}
\pswedge[linestyle = none, fillcolor = black, fillstyle = solid]
(0.000000,0.000000){1.000000}{45.000000}{60.000000}
\pswedge[linestyle = none, fillcolor = white, fillstyle = solid]
(0.000000,0.000000){1.000000}{60.000000}{75.000000}
\pswedge[linestyle = none, fillcolor = black, fillstyle = solid]
(0.000000,0.000000){1.000000}{75.000000}{90.000000}
\pswedge[linestyle = none, fillcolor = white, fillstyle = solid]
(0.000000,0.000000){1.000000}{90.000000}{105.000000}
\pswedge[linestyle = none, fillcolor = black, fillstyle = solid]
(0.000000,0.000000){1.000000}{105.000000}{120.000000}
\pswedge[linestyle = none, fillcolor = white, fillstyle = solid]
(0.000000,0.000000){1.000000}{120.000000}{135.000000}
\pswedge[linestyle = none, fillcolor = black, fillstyle = solid]
(0.000000,0.000000){1.000000}{135.000000}{150.000000}
\pswedge[linestyle = none, fillcolor = white, fillstyle = solid]
(0.000000,0.000000){1.000000}{150.000000}{165.000000}
\pswedge[linestyle = none, fillcolor = black, fillstyle = solid]
(0.000000,0.000000){1.000000}{165.000000}{180.000000}
\pswedge[linestyle = none, fillcolor = white, fillstyle = solid]
(0.000000,0.000000){1.000000}{180.000000}{195.000000}
\pswedge[linestyle = none, fillcolor = black, fillstyle = solid]
(0.000000,0.000000){1.000000}{195.000000}{210.000000}
\pswedge[linestyle = none, fillcolor = white, fillstyle = solid]
(0.000000,0.000000){1.000000}{210.000000}{225.000000}
\pswedge[linestyle = none, fillcolor = black, fillstyle = solid]
(0.000000,0.000000){1.000000}{225.000000}{240.000000}
\pswedge[linestyle = none, fillcolor = white, fillstyle = solid]
(0.000000,0.000000){1.000000}{240.000000}{255.000000}
\pswedge[linestyle = none, fillcolor = black, fillstyle = solid]
(0.000000,0.000000){1.000000}{255.000000}{270.000000}
\pswedge[linestyle = none, fillcolor = white, fillstyle = solid]
(0.000000,0.000000){1.000000}{270.000000}{285.000000}
\pswedge[linestyle = none, fillcolor = black, fillstyle = solid]
(0.000000,0.000000){1.000000}{285.000000}{300.000000}
\pswedge[linestyle = none, fillcolor = white, fillstyle = solid]
(0.000000,0.000000){1.000000}{300.000000}{315.000000}
\pswedge[linestyle = none, fillcolor = black, fillstyle = solid]
(0.000000,0.000000){1.000000}{315.000000}{330.000000}
\pswedge[linestyle = none, fillcolor = white, fillstyle = solid]
(0.000000,0.000000){1.000000}{330.000000}{345.000000}
\pswedge[linestyle = none, fillcolor = black, fillstyle = solid]
(0.000000,0.000000){1.000000}{345.000000}{360.000000}
\endpspicture
}

\def\FigureArchimedesTwo{%
\psset{unit=1cm}
\pspicture(-2,-0.5)(5,1)
\pswedge[linestyle = solid, linearc = 0.1cm,cornersize = absolute,linecolor=black,fillcolor = black, fillstyle = solid]
(0.000000,0.000000){1.000000}{82.500000}{97.500000}
\pswedge[linestyle = solid, linearc = 0.1cm,cornersize = absolute,linecolor=black,fillcolor = white, fillstyle = solid]
(0.130160,1.000000){1.000000}{262.500000}{277.500000}
\pswedge[linestyle = solid, linearc = 0.1cm,cornersize = absolute,linecolor=black,fillcolor = black, fillstyle = solid]
(0.260319,0.000000){1.000000}{82.500000}{97.500000}
\pswedge[linestyle = solid, linearc = 0.1cm,cornersize = absolute,linecolor=black,fillcolor = white, fillstyle = solid]
(0.390479,1.000000){1.000000}{262.500000}{277.500000}
\pswedge[linestyle = solid, linearc = 0.1cm,cornersize = absolute,linecolor=black,fillcolor = black, fillstyle = solid]
(0.520639,0.000000){1.000000}{82.500000}{97.500000}
\pswedge[linestyle = solid, linearc = 0.1cm,cornersize = absolute,linecolor=black,fillcolor = white, fillstyle = solid]
(0.650798,1.000000){1.000000}{262.500000}{277.500000}
\pswedge[linestyle = solid, linearc = 0.1cm,cornersize = absolute,linecolor=black,fillcolor = black, fillstyle = solid]
(0.780958,0.000000){1.000000}{82.500000}{97.500000}
\pswedge[linestyle = solid, linearc = 0.1cm,cornersize = absolute,linecolor=black,fillcolor = white, fillstyle = solid]
(0.911118,1.000000){1.000000}{262.500000}{277.500000}
\pswedge[linestyle = solid, linearc = 0.1cm,cornersize = absolute,linecolor=black,fillcolor = black, fillstyle = solid]
(1.041277,0.000000){1.000000}{82.500000}{97.500000}
\pswedge[linestyle = solid, linearc = 0.1cm,cornersize = absolute,linecolor=black,fillcolor = white, fillstyle = solid]
(1.171437,1.000000){1.000000}{262.500000}{277.500000}
\pswedge[linestyle = solid, linearc = 0.1cm,cornersize = absolute,linecolor=black,fillcolor = black, fillstyle = solid]
(1.301596,0.000000){1.000000}{82.500000}{97.500000}
\pswedge[linestyle = solid, linearc = 0.1cm,cornersize = absolute,linecolor=black,fillcolor = white, fillstyle = solid]
(1.431756,1.000000){1.000000}{262.500000}{277.500000}
\pswedge[linestyle = solid, linearc = 0.1cm,cornersize = absolute,linecolor=black,fillcolor = black, fillstyle = solid]
(1.561916,0.000000){1.000000}{82.500000}{97.500000}
\pswedge[linestyle = solid, linearc = 0.1cm,cornersize = absolute,linecolor=black,fillcolor = white, fillstyle = solid]
(1.692075,1.000000){1.000000}{262.500000}{277.500000}
\pswedge[linestyle = solid, linearc = 0.1cm,cornersize = absolute,linecolor=black,fillcolor = black, fillstyle = solid]
(1.822235,0.000000){1.000000}{82.500000}{97.500000}
\pswedge[linestyle = solid, linearc = 0.1cm,cornersize = absolute,linecolor=black,fillcolor = white, fillstyle = solid]
(1.952395,1.000000){1.000000}{262.500000}{277.500000}
\pswedge[linestyle = solid, linearc = 0.1cm,cornersize = absolute,linecolor=black,fillcolor = black, fillstyle = solid]
(2.082554,0.000000){1.000000}{82.500000}{97.500000}
\pswedge[linestyle = solid, linearc = 0.1cm,cornersize = absolute,linecolor=black,fillcolor = white, fillstyle = solid]
(2.212714,1.000000){1.000000}{262.500000}{277.500000}
\pswedge[linestyle = solid, linearc = 0.1cm,cornersize = absolute,linecolor=black,fillcolor = black, fillstyle = solid]
(2.342874,0.000000){1.000000}{82.500000}{97.500000}
\pswedge[linestyle = solid, linearc = 0.1cm,cornersize = absolute,linecolor=black,fillcolor = white, fillstyle = solid]
(2.473033,1.000000){1.000000}{262.500000}{277.500000}
\pswedge[linestyle = solid, linearc = 0.1cm,cornersize = absolute,linecolor=black,fillcolor = black, fillstyle = solid]
(2.603193,0.000000){1.000000}{82.500000}{97.500000}
\pswedge[linestyle = solid, linearc = 0.1cm,cornersize = absolute,linecolor=black,fillcolor = white, fillstyle = solid]
(2.733353,1.000000){1.000000}{262.500000}{277.500000}
\pswedge[linestyle = solid, linearc = 0.1cm,cornersize = absolute,linecolor=black,fillcolor = black, fillstyle = solid]
(2.863512,0.000000){1.000000}{82.500000}{97.500000}
\pswedge[linestyle = solid, linearc = 0.1cm,cornersize = absolute,linecolor=black,fillcolor = white, fillstyle = solid]
(2.993672,1.000000){1.000000}{262.500000}{277.500000}
\pspolygon[fillcolor=white,linestyle=none,fillstyle=solid]
(-0.1,0)(4.0,0)(4.0,-0.1)(-0.1,-0.1)
\pspolygon[fillcolor=white,linestyle=none,fillstyle=solid]
(0.0,1)(4.0,1)(4.0,1.1)(0.0,1.1)
\endpspicture
}

\def\FigureFLMKBigRectangleB{%
\psset{unit=1cm}
\pspicture(-1.3,-0.25)(2,3.2)
\pspolygon(0,0)(2,0)(2,3)(0,3)
\pspolygon[fillcolor=lightblue,fillstyle=solid](0.8,0)(0.8,3)(2,1.5)(2,1.5)
\pspolygon[fillcolor=lightyellow,fillstyle=solid](0.8,0) (0,1.5)(0.8,3)
\psdot(0,0)
\put(-0.4,-0.35){$K$}
\psdot(2,0)
\put(2.1,-0.35){$M$}
\psdot(0,1.5)
\put(-0.4,1.4){$F$}
\psdot(2,1.5)
\put(2.1,1.4){$L$}
\psdot(0.8,0)
\put(0.6,-0.35){$H$}
\psdot(0.8,1.5)
\put(0.45,1.65){$G$}
\psdot(0,3)
\put(-0.4,2.9){$k$}
\psdot(0.8,3)
\put(0.6,3.15){$h$}
\psdot(2,3)
\put(2.1,2.9){$m$}
\qline(0,1.5)(2,1.5)
\endpspicture
}

\def\FigureCyclicQuadrilateral{%
\psset{unit=1cm}
\pspicture(-3,-3.2)(3,3.0)
\pscircle(0,0){3}
\pspolygon[fillstyle=solid,fillcolor=lightyellow,linestyle=none]
(-0.5,1.4)(-2.65,1.4)(-0.5,-2.95)
\pspolygon[fillstyle=solid,fillcolor=lightyellow,linestyle=none](-0.5,1.4)(2.65,1.4)(-0.5,2.95)
\psline[linecolor=red](-0.5,-2.95)(-2.65,1.4)  
\psline[linecolor=red](2.65,1.4)(-0.5,2.95) 
\qline(-2.65,1.4)(2.65,1.4) 
\qline(-0.5,-2.95)(-0.5,2.95)
\psdot(-0.5,1.4)
\put(-0.9,1.1){$O$}
\psdot(2.65,1.4)
\put(2.8,1.35){$B$}
\psdot(-0.5,-2.95)
\put(-0.7,-3.3){$A$}
\put(-0.7,3.0){$C$}
\psdot(-2.65,1.4)
\put(-3.1,1.3){$D$}
\qline(-0.5,-2.95)(2.65,1.4) 
\qline(-0.5,2.95)(-2.65,1.4)  
\endpspicture
}

\def\FigureKupfferOne{%
\psset{unit=1cm}
\pspicture(-3,-3.2)(3,3.2)
\pscircle(0,0){3}
\pspolygon[fillstyle=solid,fillcolor=lightyellow,linestyle=none]
(-0.5,1.4)(-2.65,1.4)(-0.5,-2.95)
\pspolygon[fillstyle=solid,fillcolor=lightyellow,linestyle=none](-0.5,1.4)(2.65,1.4)(-0.5,2.95)
\psline[linecolor=red](-0.5,-2.95)(-2.65,1.4)  
\psline[linecolor=red](2.65,1.4)(-0.5,2.95) 
\qline(-2.65,1.4)(2.65,1.4) 
\qline(-0.5,-2.95)(-0.5,2.95)
\psdot(-0.5,1.4)
\put(-0.9,1.1){$O$}
\psdot(2.65,1.4)
\put(2.8,1.35){$B$}
\psdot(-0.5,-2.95)
\put(-0.7,-3.3){$A$}
\psdot(-0.5,-0.8)
\put(-0.8,-0.9){$a$}
\psdot(0.8,1.4)
\put(0.8,1.55){$b$}
\psdot(-0.5,2.95)
\put(-0.7,3.1){$\hat b$}
\psdot(-2.65,1.4)
\put(-3.1,1.3){$\hat a$}
\qline(-0.5,-0.8)(0.8,1.4) 
\qline(-0.5,-2.95)(2.65,1.4) 
\qline(-0.5,2.95)(-2.65,1.4)  
\endpspicture
}

\def\FigureKupfferTwo{%
\psset{unit=0.8cm}
\pspicture(-1.3,-8)(8,6)
\qline(0,0)(8,0)
\qline(0,-8)(0,5.7)
\psdot(0,-7.5) 
\put(-0.6,-7.6){$D^\prime$}
\psdot(0,0)  
\put(-0.6,-0.1){$O$}
\psdot(0,3)  
\put(-0.6,2.9){$B^\prime$}
\psdot(0,2)  
\put(-0.6,1.9){$A^\prime$}
\psdot(3,0)  
\put(2.5,-0.4){$C$}
\psdot(0,5) 
\put(-0.6,4.9){$C^\prime$}
\psdot(5,0) 
\put(5.55,-0.4){$B$}
\psdot(7.5,0) 
\put(7.6,-0.4){$A$}
\psline[linestyle=solid,linecolor=red](0,2)(5,0) 
\psline[linestyle=solid,linecolor=red](0,3)(7.5,0) 
\psline[linestyle=solid,linecolor=blue](0,2)(3,0) 
\psline[linestyle=dashed,linecolor=blue](3,0)(4.44,-0.92) 
\psline[linestyle=solid,linecolor=blue](0,5)(7.5,0) 
\psline[linestyle=solid,linecolor=green](0,3)(3,0) 
\psline[linestyle=dashed,linecolor=green](3,0)(5.2,-2.2) 
\psline[linestyle=solid,linecolor=green](0,5)(5,0) 
\psline[linestyle=dashed,linecolor=green](5,0)(6.3,-1.3) 
\psline[linestyle=dotted,linecolor =gray](0,-7.5)(5,0)  
\psline[linestyle=dotted,linecolor =gray](5,0)(5.76,1.17)  
\psline[linestyle=dotted,linecolor =gray](0,-7.5)(3,0)  
\psline[linestyle=dotted,linecolor =gray](3,0)(3.6,1.6)  
\psline[linestyle=dotted,linecolor =gray](0,-7.5)(7.5,0)  
\endpspicture
}

\def\FigureProportionOne{%
\psset{unit=0.8cm}
\pspicture(-1.3,-0.5)(8,6)
\pspolygon[fillstyle=solid,fillcolor=lightyellow](0,0)(3.48,3.48)(5,0)
\pspolygon[fillstyle=solid,fillcolor=yellow](0,0)(2.13,2.13)(3,0)
\qline(0,0)(8,0)
\qline(0,0)(4,4)
\qline(0,0)(0,5.7)
\psdot(0,0)  
\put(-0.6,-0.1){$A$}
\psdot(0,3)  
\put(-0.6,2.9){$d$}
\psdot(3,0)  
\put(2.9,-0.4){$c$}
\psdot(0,5) 
\put(-0.6,4.9){$D$}
\psdot(5,0) 
\put(4.9,-0.45){$C$}
\psdot(2.13,2.13) 
\put(2.35,2){$b$}
\psdot(3.48,3.48) 
\put(3.6,3.3){$B$}
\psline[linestyle=solid,linecolor=red](0,5)(3.48,3.48) 
\psline[linestyle=solid,linecolor=red](0,3)(2.13,2.13) 
\psline[linestyle=solid,linecolor=blue](0,3)(3,0) 
\psline[linestyle=solid,linecolor=blue](0,5)(5,0) 
\psline[linestyle=solid,linecolor=green](3.48,3.48)(5,0) 
\psline[linestyle=solid,linecolor=green](2.13,2.13)(3,0) 
\endpspicture
}

\def\FigureProportionDefn{%
\psset{unit=0.8cm}
\pspicture(-1.3,-0.5)(8,6)
\qline(0,0)(8,0)
\qline(0,0)(0,5.7)
\psdot(0,0)  
\put(-0.6,-0.1){$A$}
\psdot(0,3)  
\put(-0.6,2.9){$b$}
\psdot(4.5,0)  
\put(4.4,-0.4){$c$}
\psdot(0,5) 
\put(-0.6,4.9){$B$}
\psdot(7.5,0) 
\put(7.4,-0.45){$C$}
\psline[linestyle=solid,linecolor=blue](0,3)(4.5,0) 
\psline[linestyle=solid,linecolor=blue](0,5)(7.5,0) 
\endpspicture
}

\def\FigureProportionalSimilar{%
\psset{unit=0.8cm}
\pspicture(-1.3,-0.5)(5,2.7)
\qline(0,0)(4.5,0)
\psdot(0,0)  
\put(-0.6,-0.1){$A$}
\psdot(4.5,0)  
\put(4.4,-0.45){$B$}
\psdot(2.7,0) 
\put(2.7,-0.45){$b$}
\psdot(4.5,2.18) 
\put(4.6,2.2){$C$}
\qline(4.5,2.18)(4.5,0) 
\qline(0,0)(4.5,2.18) 
\psdot(2.7,1.7) 
\put(2.9,1.9){$c$}
\qline(0,0)(2.7,1.7)  
\qline(2.7,1.7)(2.7,0) 
\psdot(2.7, 1.3) 
\put(2.9,1.0){$P$}
\endpspicture
}

\def\FigureFLMKParallelogram{%
\psset{unit=1cm}
\pspicture(-1.3,-0.25)(2,3)
\pspolygon[fillcolor=lightblue,fillstyle=solid](0.8,0)(1.3,1.5)(2.5,1.5)(2,0)
\pspolygon[fillcolor=lightyellow,fillstyle=solid](0.8,0)(0,0)(0.5,1.5)(1.3,1.5)
\psdot(0,0)
\put(-0.4,-0.35){$K$}
\psdot(2,0)
\put(2.1,-0.35){$M$}
\psdot(0.5,1.5)
\put(0.1,1.4){$F$}
\psdot(2.5,1.5)
\put(2.6,1.4){$L$}
\psdot(0.8,0)
\put(0.6,-0.35){$H$}
\psdot(1.3,1.5)
\put(1.1,1.65){$G$}
\endpspicture
}

\def\FigurePasteFiveHelperG{%
\psset{unit=1cm}
\pspicture(-0.3,-0.25)(2,3.2)
\pspolygon[fillcolor=lightyellow,fillstyle=solid](0,1)(0.8,3)(0.8,0)
\pspolygon[fillcolor=lightblue,fillstyle=solid](0.8,3)(2,2)(0.8,0)
\psdot(0,1)
\psdot(0.8,3)
\psdot(2,2)
\psdot(0.8,0)
\psline[linestyle=dashed] (0.8,0)(0.8,3)
\put(-0.4,0.9){$A$}
\put(0.6,3.13){$B$}
\put(2.1,1.9){$C$}
\put(0.6,-0.35){$D$}
\endpspicture
}

\def\FigurePasteFiveHelperA{%
\psset{unit=1cm}
\pspicture(-0.3,-0.25)(2,3)
\pspolygon[linecolor=red]
(0,0)(0,3)(2,3)(2,0)
\pspolygon[fillcolor=lightblue,fillstyle=solid](0,1)(0.8,3)(2,2)(0.8,0)
\psdot(0,1)
\psdot(0.8,3)
\psdot(2,2)
\psdot(0.8,0)
\psline[linestyle=dashed] (0.8,0)(0.8,3)
\put(-0.4,0.9){$A$}
\put(0.6,3.13){$B$}
\put(2.1,1.9){$C$}
\put(0.6,-0.35){$D$}
\endpspicture
}

\def\FigurePasteFiveHelperB{%
\psset{unit=1cm}
\pspicture(-1.3,-0.25)(2,3)
\pspolygon[linecolor=red](0,0)(0,3)(2,3)(2,0)
\pspolygon[fillcolor=lightblue,fillstyle=solid](0,0)(0,1.5)(2,1.5)(2,0)
\psdot(0,0)
\put(-0.4,-0.35){$S$}
\psdot(2,0)
\put(2.1,-0.35){$R$}
\psdot(0,3)
\put(-0.4,3.13){$P$}
\psdot(2,3)
\put(2.1,3.13){$Q$}
\psdot(0,1.5)
\put(-0.4,1.4){$G$}
\psdot(2,1.5)
\put(2.1,1.4){$H$}
\endpspicture
}

\def\FigurePasteFiveHelperC{%
\psset{unit=1cm}
\pspicture(-0.3,-0.25)(2,3)
\pspolygon[linecolor=red]
(0,0)(0,3)(2,3)(2,0)
\pspolygon[fillcolor=lightblue,fillstyle=solid](0,1.5)(0.8,3)(2,1.5)(0.8,0)
\psdot(0.8,3)
\psdot(0.8,0)
\psline[linestyle=dashed] (0.8,0)(0.8,3)
\psdot(0,1.5)
\put(-0.4,1.4){$A$}
\psdot(2,1.5)
\put(2.1,1.4){$C$}
\put(0.6,3.13){$B$}
\put(0.6,-0.35){$D$}
\endpspicture
}

\def\FigurePasteFiveHelperD{%
\psset{unit=1cm}
\pspicture(-1.3,-0.25)(2,3)
\pspolygon[linecolor=red]
(0,0)(0,3)(2,3)(2,0)
\pspolygon[fillcolor=lightblue,fillstyle=solid](0,1.5)(0.8,3)(2,1.5)(0.8,0)
\psdot(0.8,3)
\psdot(0.8,0)
\psline[linestyle=dashed] (0,1.5)(2,1.5)
\psdot(0,1.5)
\put(-0.4,1.4){$A$}
\psdot(2,1.5)
\put(2.1,1.4){$C$}
\put(0.6,3.13){$B$}
\put(0.6,-0.35){$D$}
\endpspicture
}

\def\FigurePasteFiveHelperE{%
\psset{unit=1cm}
\pspicture(-1.3,-0.25)(2,3)
\pspolygon[linecolor=red]
(0,0)(0,3)(2,3)(2,0)
\pspolygon[fillcolor=lightblue,fillstyle=solid](0,1.5)(0,3)(2,1.5)(2,0)
\psdot(0,1.5)
\put(-0.4,1.4){$A$}
\psdot(2,1.5)
\put(2.1,1.4){$C$}
\psdot(0,3)
\put(-0.4,3.13){$B$}
\psdot(2,0)
\put(2.1,-0.35){$D$}
\psline[linestyle=dashed] (0,1.5)(2,1.5) 
\endpspicture
}
 
\def\FigurePasteFiveHelperF{%
\psset{unit=1cm}
\pspicture(-1.3,-0.25)(2,3)
\pspolygon[linecolor=red](0,0)(0,3)(2,3)(2,0)
\pspolygon[fillcolor=lightblue,fillstyle=solid](0,0)(0,1.5)(2,1.5)(2,0)
\psdot(0,0)
\put(-0.4,-0.35){$S$}
\psdot(2,0)
\put(2.1,-0.35){$R$}
\psdot(0,3)
\put(-0.4,3.13){$P$}
\psdot(2,3)
\put(2.1,3.13){$Q$}
\psdot(0,1.5)
\put(-0.4,1.4){$G$}
\psdot(2,1.5)
\put(2.1,1.4){$H$}
\psline[linestyle=dashed] (0,1.5)(2,0) 
\endpspicture
}

\def\FigureERFourA{%
\psset{unit=1.2cm}
\pspicture(-0.3,-0.15)(5,3.15)
\pspolygon[fillcolor=lightblue,fillstyle=solid,linestyle=none]
(0,1.32)(1.2,1.32)(1.2,3)(0,3)
\put(0.5,2.5){{\tiny b}}
\pspolygon[fillcolor=lightblue,fillstyle=solid,linestyle=none]
(2,0)(4.5,0)(4.5,0.8)(2,0.8)
\put(2.5,0.4){{\tiny b}}
\pspolygon[fillcolor=pink,fillstyle=solid,linestyle=none]
(1.2,1.32)(2,1.32)(2,3)(1.2,3)
\put(1.5,2.5){{\tiny p}}
\pspolygon[fillcolor=pink,fillstyle=solid,linestyle=none]
(2,1.32)(4.5,1.32)(4.5,0.8)(2,0.8)
\put(2.5,1.05){{\tiny p}}
\pspolygon[fillcolor=lightyellow,fillstyle=solid,linestyle=none]
(1.2,0.8)(2,0.8)(2,1.32)(1.2,1.32)
\put(1.4,1.1){{\tiny y}}
\pspolygon[fillcolor=lightgreen,fillstyle=solid,linestyle=none]
(0,0)(1.2,0)(1.2,0.8)(0,0.8)
\put(0.5,0.5){{\tiny g}}
\pspolygon[fillcolor=red,fillstyle=solid,linestyle=none]
(0,0.8)(1.2,0.8)(1.2,1.32)(0,1.32)
\put(0.5,1.05){{\tiny r}}
\pspolygon[fillcolor=red,fillstyle=solid,linestyle=none]
(1.2,0)(2,0)(2,0.8)(1.2,0.8)
\put(1.5,0.4){{\tiny r}}
\qline(0,0)(0,3)
\qline(0,0)(4.5,3)
\qline(0,0)(4.5,0)
\qline(0,3)(4.5,3)
\qline(4.5,0)(4.5,3)
\qline(2,0)(2,3)
\qline(0,0.8)(4.5,0.8)
\qline(0,1.32)(4.5,1.32)
\qline(1.2,0)(1.2,3)
\psdot(0,3)
\psdot(0,1.32)
\psdot(0,0)
\psdot(1.2,3)
\psdot(2,3)
\psdot(4.5,3)
\psdot(4.5,1.32)
\psdot(4.5,0.8)
\psdot(4.5,0)
\psdot(2,0)
\psdot(1.2,0.8)
\psdot(2,0.8)
\psdot(1.2,1.32)
\psdot(2,1.32)
\psdot(0,0.8)
\psdot(1.2,0)
\endpspicture
}

\def\FigureERFourB{%
\psset{unit=1.2cm}
\pspicture(-0.3,-0.15)(5,3.15)
\pspolygon[fillcolor=pink,fillstyle=solid,linestyle=none]
(0,1.32)(0.8,1.32)(0.8,3)(0,3)
\put(0.35,2.5){{\tiny p}}
\pspolygon[fillcolor=pink,fillstyle=solid,linestyle=none]
(2,0)(4.5,0)(4.5,0.52)(2,0.52)
\put(2.5, 0.2){{\tiny p}}
\pspolygon[fillcolor=lightblue,fillstyle=solid,linestyle=none]
(0.8,1.32)(2,1.32)(2,3)(0.8,3)
\put(1.35,2.5){{\tiny b}}
\pspolygon[fillcolor=lightblue,fillstyle=solid,linestyle=none]
(2,1.32)(4.5,1.32)(4.5,0.52)(2,0.52)
\put(2.5, 0.9){{\tiny b}}
\pspolygon[fillcolor=lightgreen,fillstyle=solid,linestyle=none]
(0.8,0.52)(2,0.52)(2,1.32)(0.8,1.32)
\put(1.05,0.9){{\tiny g}}
\pspolygon[fillcolor=lightyellow,fillstyle=solid,linestyle=none]
(0,0)(0.8,0)(0.8,0.52)(0,0.52)
\put(0.15,0.3){{\tiny y}}
\pspolygon[fillcolor=red,fillstyle=solid,linestyle=none]
(0,0.52)(0.8,0.52)(0.8,1.32)(0,1.32)
\put(0.35,0.9){{\tiny r}}
\pspolygon[fillcolor=red,fillstyle=solid,linestyle=none]
(0.8,0)(2,0)(2,0.52)(0.8,0.52)
\put(1.35, 0.2){{\tiny r}}
\qline(0,0)(0,3)
\qline(0,0)(4.5,3)
\qline(0,0)(4.5,0)
\qline(0,3)(4.5,3)
\qline(4.5,0)(4.5,3)
\qline(2,0)(2,3)
\qline(0,0.52)(4.5,0.52)
\qline(0,1.32)(4.5,1.32)
\qline(0.8,0)(0.8,3)
\psdot(0,3)
\psdot(0,1.32)
\psdot(0,0)
\psdot(0.8,3)
\psdot(2,3)
\psdot(4.5,3)
\psdot(4.5,1.32)
\psdot(4.5,0.52)
\psdot(4.5,0)
\psdot(2,0)
\psdot(0.8,0.52)
\psdot(2,0.52)
\psdot(0.8,1.32)
\psdot(2,1.32)
\psdot(0.8,0)
\psdot(0,0.52)
\endpspicture
}

\def\FigureOneThirtyFiveWithDiagonals{%
\psset{unit=1.3cm}
\pspicture(0,-0.5)(4,2)
\psdot(0.5,0)
\put(0.4,-0.3){$B$}
\psdot(0,2)
\put(-0.1,2.15){$A$}
\psdot(2,2)
\put(1.9,2.15){$D$}
\psdot(2.5,2)
\put(2.3,2.15){$E$}
\psdot(4.5,2)
\put(4.3,2.15){$F$}
\psdot(2.5,0)
\put(2.3,-0.3){$C$}
\psline[linecolor=red](0,2)(2.5,0)  
\psline[linecolor=red](0.5,0)(4.5,2) 
\qline(0,2)(4.5,2)  
\qline(0,2)(0.5,0)  
\qline(0.5,0)(2.5,2)  
\qline(2,2)(2.5,0)  
\qline(4.5,2)(2.5,0) 
\qline(0.5,0)(2.5,0) 
\psdot(2.1,1.6)
\put(2.22,1.5){$G$}
\endpspicture
}

\def\FigureOneThirtyFiveColored{%
\psset{unit=1.5cm}
\pspicture(0,-0.5)(4,2.5)
\pspolygon[fillstyle=solid,fillcolor=yellow,linestyle=none]
(0,2)(0.5,0)(2.5,2)
\pspolygon[fillstyle=solid,fillcolor=yellow,linestyle=none]
(2,2)(4.5,2)(2.5,0)
\pspolygon[fillstyle=solid,fillcolor=red,linestyle=none]
(2,2)(2.5,2)(2.1,1.6)
\psdot(0.5,0)
\put(0.4,-0.3){$B$}
\psdot(0,2)
\put(-0.1,2.15){$A$}
\psdot(2,2)
\put(1.9,2.15){$D$}
\psdot(2.5,2)
\put(2.3,2.15){$E$}
\psdot(4.5,2)
\put(4.3,2.15){$F$}
\psdot(2.5,0)
\put(2.3,-0.3){$C$}
\qline(0,2)(4.5,2)  
\qline(0,2)(0.5,0)  
\qline(0.5,0)(2.5,2)  
\qline(2,2)(2.5,0)  
\qline(4.5,2)(2.5,0) 
\qline(0.5,0)(2.5,0) 
\psdot(2.1,1.6)
\put(2.22,1.5){$G$}
\endpspicture
}

\def\FigureOneFortyThree{%
\psset{unit=0.8cm}
\pspicture(0,-0.5)(7,4.6)
\pspolygon[fillcolor=yellow,fillstyle=solid]
(0,0)(0.6,3)(2.6,3)(2,0)
\pspolygon[fillcolor=yellow,fillstyle=solid]
(2.6,3)(6.98,3)(7.26,4.4)(2.88,4.4)
\pspolygon
(0,0)(0.88,4.4)(7.26,4.4)(6.38,0)
\qline(6.38,0)(0.88,4.4)  
\psdot(0,0)
\put(-0.1,-0.5){$B$}
\psdot(2,0)
\put(1.8,-0.5){$G$}
\psdot(6.38,0)
\put(6.3,-0.5){$C$}
\psdot(0.88,4.4)
\put(0.75,4.6){$A$}
\psdot(2.88,4.4)
\put(2.75,4.6){$H$}
\psdot(7.26,4.4)
\put(7.13,4.6){$D$}
\psdot(0.6,3)
\put(0.1,2.9){$E$}
\psdot(6.98,3)
\put(7.2,2.9){$F$}
\psdot(2.6,3)
\put(2.8,3.2){$K$}
\endpspicture
}

\def\FigureHalvesOfRectangles{%
\psset{unit=1.5cm}
\pspicture(0,-0.15)(6,4.15)
\pspolygon[fillcolor=pink,fillstyle=solid]
(0,4)(2,4)(2,1.32)(0,1.32)
\pspolygon[fillcolor=lightblue,fillstyle=solid]
(2,1.32)(2,0)(6,0)(6,1.32)
\pspolygon[fillcolor=lightyellow,fillstyle=solid]
(2,1.32)(1,1.32)(1,0.66)(2,0.66)
\qline(0,0)(0,4)
\qline(0,0)(6,4)
\qline(0,0)(6,0)
\qline(0,4)(6,4)
\qline(6,0)(6,4)
\qline(2,0)(2,4)
\qline(1,0.66)(6,0.66)
\qline(0,1.32)(6,1.32)
\qline(1,0.66)(1,4)
\psdot(0,4)
\put(-0.3,3.9){$A$}
\psdot(0,1.32)
\put(-0.3,1.22){$D$}
\psdot(0,0)
\put(-0.3,-0.1){$E$}
\psdot(1,4)
\put(0.9,4.1){$P$}
\psdot(2,4)
\put(1.9,4.1){$B$}
\psdot(6,4)
\put(6.1,3.9){$F$}
\psdot(6,1.32)
\put(6.1,1.22){$b$}
\psdot(6,0.66)
\put(6.1,0.56){$q$}
\psdot(6,0)
\put(6.1,-0.1){$c$}
\psdot(2,0)
\put(1.9,-0.25){$d$}
\psdot(1,0.66)
\put(0.75,0.72){$R$}
\psdot(2,0.66)
\put(2.1,0.75){$p$}
\psdot(1,1.32)
\put(0.75,1.42){$Q$}
\psdot(2,1.32)
\put(1.75,1.42){$C$}
\endpspicture
}

\def\FigureEqualRQ{%
\psset{unit=1.5cm}
\pspicture(0,-0.15)(6,4.5)
\pspolygon[linecolor=red,fillcolor=pink,fillstyle=solid]
(0,1.75)(1,4)(2,3.57)(1,1.32)
\pspolygon[linecolor=blue,fillcolor=lightblue,fillstyle=solid]
(2,0.66)(2.14,1.32)(6,0.66)(5.86,0)
\qline(0,0)(0,4)
\qline(0,0)(6,4)
\qline(0,0)(6,0)
\qline(0,4)(6,4)
\qline(6,0)(6,4)
\qline(2,0)(2,4)
\qline(2,0.66)(6,0.66)
\qline(0,1.32)(6,1.32)
\qline(1,1.32)(1,4)
\psdot(0,1.75)
\put(-0.3,1.7){$A$}
\psdot(2,3.57)
\put(2.1,3.5){$C$}
\psdot(5.86,0)
\put(5.8,-0.25){$b$}
\psdot(2.2,1.32)
\put(2.2,1.05){$d$}
\psdot(0,4)
\psdot(0,1.32)
\psdot(0,0)
\put(-0.3,-0.1){$E$}
\psdot(1,4)
\put(0.9,4.1){$B$}
\psdot(2,4)
\psdot(6,4)
\put(6.1,3.9){$F$}
\psdot(6,1.32)
\psdot(6,0.66)
\put(6.1,0.56){$a$}
\psdot(6,0)
\psdot(2,0)
\psdot(2,0.66)
\put(1.77,0.6){$c$}
\psdot(1,1.32)
\put(0.9,1.05){$D$}
\psdot(2,1.32)
\put(1.75,1.42){$H$}
\endpspicture
}

\def\FigureHalvesOfEquals{%
\psset{unit=1.5cm}
\pspicture(0,-0.15)(6,4.15)
\pspolygon[linecolor=red,fillcolor=pink,fillstyle=solid]
(0,3.5)(1,4)(2,3)(1,1.32)
\pspolygon[linecolor=blue,fillstyle=solid,fillcolor=lightblue]
(2,0.66)(3,1.32)(6,0.66)(2.5,0)
\qline(0,0)(0,4)
\qline(0,0)(6,4)
\qline(0,0)(6,0)
\qline(0,4)(6,4)
\qline(6,0)(6,4)
\qline(2,0)(2,4)
\qline(1,0.66)(6,0.66)
\qline(0,1.32)(6,1.32)
\qline(1,0.66)(1,4)
\psdot(0,4)
\put(-0.3,3.9){$A$}
\psdot(0,1.32)
\put(-0.3,1.22){$D$}
\psdot(0,0)
\put(-0.3,-0.1){$E$}
\psdot(1,4)
\put(0.9,4.1){$P$}
\psdot(2,4)
\put(1.9,4.1){$B$}
\psdot(6,4)
\put(6.1,3.9){$F$}
\psdot(6,1.32)
\put(6.1,1.22){$b$}
\psdot(6,0.66)
\put(6.1,0.56){$q$}
\psdot(6,0)
\put(6.1,-0.1){$c$}
\psdot(2,0)
\put(1.9,-0.25){$d$}
\psdot(1,0.66)
\put(0.75,0.72){$R$}
\psdot(2,0.66)
\put(1.8,0.75){$p$}
\psdot(1,1.32)
\put(0.75,1.12){$Q$}
\psdot(2,1.32)
\put(1.75,1.42){$C$}
\endpspicture
}

\def\EqualRectanglesFigureOne{%
\pspicture(-0.15,-0.15)(1.7,1.5)
\pspolygon[fillcolor=yellow,fillstyle=solid,linestyle=none]
(0,1.5)(0.6,1.5)(0.6,0.3)(0.0,0.3)
\pspolygon[fillcolor=yellow,fillstyle=solid,linestyle=none]
(0.6,0.3)(0.6,0)(1.5,0)(1.5,0.3)
\qline(0.0,0.3)(1.5,0.3)  
\qline(0.6,0)(1.5,0)      
\qline(0.0,1.5)(0.6,1.5)  
\qline(0.0,0.3)(0.0,1.5)  
\qline(0.6,0)(0.6,1.5)    
\qline(1.5,0)(1.5,0.3)    
\psdot(0.0,0.3) 
\put(-0.15,0.3){$G$}
\psdot(1.5,0.3)  
\put(1.55,0.3){$M$}
\psdot(1.5,0)  
\put(1.5,-0.12){$L$}
\psdot(0.6,0.3)  
\put(0.63,0.35){$B$}
\psdot(0.6,0) 
\put(0.63,-0.12){$A$}
\psdot(0.6,1.5) 
\put(0.63,1.55){$E$}
\psdot(0,1.5) 
\put(-0.15,1.55){$F$}
\endpspicture
}

\def\EqualRectanglesFigureTwo{%
\pspicture(-0.2,-0.15)(1.7,1.5)
\pspolygon[fillcolor=yellow,fillstyle=solid,linestyle=none]
(0,1.5)(0.6,1.5)(0.6,0.3)(0.0,0.3)
\pspolygon[fillcolor=yellow,fillstyle=solid,linestyle=none]
(0.6,0.3)(0.6,0)(1.5,0)(1.5,0.3)
\qline(0.0,0.3)(1.5,0.3)  
\qline(0.6,0)(1.5,0)      
\qline(0.0,1.5)(0.6,1.5)  
\qline(0.0,0.3)(0.0,1.5)  
\qline(0.6,0)(0.6,1.5)    
\qline(1.5,0)(1.5,0.3)    
\psdot(0.0,0.3) 
\put(-0.15,0.3){$G$}
\psdot(1.5,0.3)  
\put(1.55,0.3){$M$}
\psdot(1.5,0)  
\put(1.5,-0.12){$L$}
\psdot(0.6,0.3)  
\put(0.48,0.35){$B$}
\psdot(0.6,0) 
\put(0.63,-0.12){$A$}
\psdot(0.6,1.5) 
\put(0.63,1.55){$E$}
\psdot(0,1.5) 
\put(-0.15,1.55){$F$}
\psline[linecolor=red](0,0)(0,0.3) 
\psline[linecolor=red](0,0)(0.6,0) 
\psline[linecolor=red](0.6,1.5)(1.5,1.5) 
\psline[linecolor=red](1.5,1.5)(1.5,0.3) 
\psdot(1.5,1.5) 
\put(1.5,1.55){$K$}
\psdot(0,0) 
\put(-0.15,-0.12){$H$}
\psline[linecolor=blue,linestyle=dashed](0,0)(0.6,0.3)
\psline[linecolor=blue,linestyle=dashed](1.5,1.5)(0.6,0.3)
\endpspicture
}

\def\EqualRectanglesFigureThree{%
\pspicture(-0.2,-0.15)(1.7,0.9)
\pspolygon[fillcolor=yellow,fillstyle=solid,linestyle=none]
(0,0.75)(0.6,0.75)(0.6,0.3)(0.0,0.3)
\pspolygon[fillcolor=yellow,fillstyle=solid,linestyle=none]
(0.6,0.3)(0.6,0)(1.5,0)(1.5,0.3)
\qline(0.0,0.3)(1.5,0.3)  
\qline(0.6,0)(1.5,0)      
\qline(0.0,0.75)(0.6,0.75)  
\qline(0.0,0.3)(0.0,0.75)  
\qline(0.6,0)(0.6,0.75)    
\qline(1.5,0)(1.5,0.3)    
\psdot(0.0,0.3) 
\put(-0.15,0.3){$G$}
\psdot(1.5,0.3)  
\put(1.55,0.3){$M$}
\psdot(1.5,0)  
\put(1.5,-0.12){$L$}
\psdot(0.6,0.3)  
\put(0.48,0.35){$B$}
\psdot(0.6,0) 
\put(0.63,-0.12){$A$}
\psdot(0.6,0.75) 
\put(0.63,0.8){$E$}
\psdot(0,0.75) 
\put(-0.15,0.8){$F$}
\psline[linecolor=red](0,0)(0,0.3) 
\psline[linecolor=red](0,0)(0.6,0) 
\psline[linecolor=red](0.6,0.75)(1.5,0.75) 
\psline[linecolor=red](1.5,0.75)(1.5,0.3) 
\psdot(1.5,0.75) 
\put(1.5,0.8){$K$}
\psdot(0,0) 
\put(-0.15,-0.12){$H$}
\psline[linecolor=blue,linestyle=dashed](0,0)(1.5,0.75)
\endpspicture
}

\def\FigureABCDequalsBCDA{%
\pspicture(-0.2,-0.15)(1.7,1.0)
\pspolygon[fillcolor=yellow,fillstyle=solid,linestyle=solid]
(0,0.75)(0.75,0.75)(0.75,0.375)(0.0,0.375)
\pspolygon[fillcolor=yellow,fillstyle=solid,linestyle=solid]
(0.75,0.375)(0.75,0)(1.5,0)(1.5,0.375)
\psdot(0.0,0.375) 
\put(-0.15,0.375){$C$}
\psdot(1.5,0.375)  
\put(1.55,0.375){$C^\prime$}
\psdot(1.5,0)  
\put(1.5,-0.12){$D^\prime$}
\psdot(0.75,0.375)  
\put(0.63,0.4){$B$}
\psdot(0.75,0) 
\put(0.63,-0.12){$A^\prime$}
\psdot(0.75,0.75) 
\put(0.63,0.8){$A$}
\psdot(0,0.75) 
\put(-0.15,0.8){$D$}
\psline[linecolor=red](0,0)(0,0.375) 
\psline[linecolor=red](0,0)(0.75,0) 
\psline[linecolor=red](0.75,0.75)(1.5,0.75) 
\psline[linecolor=red](1.5,0.75)(1.5,0.375) 
\psdot(1.5,0.75) 
\put(1.5,0.8){$K$}
\psdot(0,0) 
\put(-0.15,-0.12){$H$}
\psline[linecolor=blue,linestyle=dashed](0,0)(1.5,0.75)
\endpspicture
}

\def\FigureERTransitive{%
\pspicture(-0.2,-0.15)(2.0,1.4)
\pspolygon[fillcolor=yellow,fillstyle=solid,linestyle=none]
(0,0.75)(0.6,0.75)(0.6,0.3)(0.0,0.3)
\pspolygon[fillcolor=yellow,fillstyle=solid,linestyle=none]
(0.6,0.3)(0.6,0)(1.5,0)(1.5,0.3) 
\pspolygon[fillcolor=yellow,fillstyle=solid,linestyle=solid]
(1.5,1.2) (1.8,1.2)(1.8,0.3) (1.5,0.3) 
\qline(0.0,0.3)(1.5,0.3)  
\qline(0.6,0)(1.8,0)      
\qline(0.0,0.75)(0.6,0.75)  
\qline(0.0,0.3)(0.0,0.75)  
\qline(0.6,0)(0.6,0.75)    
\psline[linecolor=red](0.6,0.75)(0.6,1.2)    
\qline(1.5,0)(1.5,0.3)    
\psdot(0.6,1.2)  
\put(0.63,1.24){$D$}
\psdot(0.0,0.3) 
\put(-0.15,0.3){$G$}
\psdot(1.5,0.3)  
\put(1.55,0.33){$M$}
\psdot(1.5,0)  
\put(1.5,-0.12){$L$}
\psdot(1.8,0) 
\put(1.8,-0.12){$P$}
\psdot(0.6,0.3)  
\put(0.48,0.35){$B$}
\psdot(0.6,0) 
\put(0.63,-0.12){$A$}
\psdot(0.6,0.75) 
\put(0.63,0.8){$E$}
\psdot(0,0.75) 
\put(-0.15,0.8){$F$}
\psline[linecolor=red](0,0)(0,0.3) 
\psline[linecolor=red](0,0)(0.6,0) 
\psline[linecolor=red](0.6,0.75)(1.5,0.75) 
\psline[linecolor=red](0.6,1.2)(1.5,1.2) 
\psline[linecolor=red](1.8,0.3)(1.8,0) 
\psdot(1.8,0.3)   
\put(1.85,0.3){$Q$}
\psdot(1.8,1.2) 
\put(1.83,1.24){$R$}
\psdot(1.5,1.2) 
\put(1.5,1.24){$C$}
\psdot(1.5,0.75) 
\put(1.52,0.8){$K$}
\psdot(0,0) 
\put(-0.15,-0.12){$H$}
\psline[linecolor=blue,linestyle=dashed](0,0)(1.5,0.75)
\psline[linecolor=blue,linestyle=dashed](1.8,0)(0.6,1.2)
\put(0.3,0.5){$1$}
\put(1.0,0.12){$2$}
\put(1.63, 0.7){$3$}
\endpspicture
}

\def\FigureETpermutation{%
\pspicture(1.2,2.4)(0,-0.15)
\psset{unit=1cm}
\pspolygon[fillstyle=solid,fillcolor=yellow,linewidth=1pt,linecolor=gray](0.00,6.25)(1.79,6.25)(1.79,2.95)(0.00,2.95)
\pspolygon[fillstyle=solid,fillcolor=yellow,linewidth=1pt,linecolor=gray](1.79,2.95)(3.79,2.95)(3.79,0.00)(1.79,0.00)
\pspolygon[fillstyle=solid,fillcolor=lightblue,linewidth=1pt,linecolor=gray](1.79,6.25)(3.79,6.25)(3.79,2.95)(1.79,2.95)
\pspolygon[fillstyle=solid,fillcolor=lightblue,linewidth=1pt,linecolor=gray](0.00,2.95)(1.79,2.95)(1.79,0.00)(0.00,0.00)
\pspolygon[fillstyle=solid,linewidth=1pt,fillcolor=lightgreen](1.79,6.25)(1.79,2.95)(0.00,3.85)
\pspolygon[fillstyle=solid,linewidth=1pt,fillcolor=lightgreen](3.26,0.00)(1.79,2.95)(3.79,2.95)
\pspolygon[fillstyle=solid,linewidth=1pt,fillcolor=pink](0.00,3.85)(1.79,2.95)(0.00,2.95)
\pspolygon[fillstyle=solid,linewidth=1pt,fillcolor=pink](1.79,2.95)(1.79,0.00)(3.26,0.00)
\psline[linecolor = blue,linestyle=dashed](0.00,0.00)(3.79,6.25)
\put(1.789854,6.353260){$A$}
\put(0.000000,6.353260){$K$}
\put(3.789854,6.353260){$G$}
\put(-0.450000,3.845684){$C$}
\put(-0.450000,2.953260){$D$}
\put(-0.450000,0.000000){$E$}
\put(2.0,3.02){$B$}
\put(1.789854,-0.300000){$F$}
\put(3.789854,-0.100000){$L$}
\put(3.262354,-0.300000){$A^\prime$}
\put(3.889854,3.003260){$C^\prime$}
\put(1.5,4.5){$c$}
\put(0.75,4.5){$b$}
\put(0.84,3.54){$a$}
\put(2.6,1.6){$c$}
\put(3.3,1.6){$b$}
\put(2.9,2.7){$a$}
\endpspicture
}

\def\FigureEqualTriangles{%
\pspicture(1.2,2.4)(0,-0.15)
\psset{unit=1cm}
\pspolygon[fillstyle=solid,fillcolor=white,linewidth=1pt,linecolor=black](0.00,6.25)(1.79,6.25)(1.79,2.95)(0.00,2.95)
\pspolygon[fillstyle=solid,fillcolor=white,linewidth=1pt,linecolor=black](1.79,2.95)(3.79,2.95)(3.79,0.00)(1.79,0.00)
\pspolygon[fillstyle=solid,fillcolor=lightblue,linewidth=1pt,linecolor=gray](1.79,6.25)(3.79,6.25)(3.79,2.95)(1.79,2.95)
\pspolygon[fillstyle=solid,fillcolor=lightblue,linewidth=1pt,linecolor=gray](0.00,2.95)(1.79,2.95)(1.79,0.00)(0.00,0.00)
\pspolygon[fillstyle=solid,linewidth=1pt,fillcolor=lightgreen](1.79,6.25)(1.79,2.95)(0.00,3.85)
\pspolygon[fillstyle=solid,linewidth=1pt,fillcolor=lightgreen](3.26,0.00)(1.79,2.95)(3.79,2.95)
\pspolygon[fillstyle=solid,linewidth=1pt,fillcolor=white](0.00,3.85)(1.79,2.95)(0.00,2.95)
\pspolygon[fillstyle=solid,linewidth=1pt,fillcolor=white](1.79,2.95)(1.79,0.00)(3.26,0.00)
\psline[linecolor = blue,linestyle=dashed](0.00,0.00)(3.79,6.25)
\put(1.789854,6.353260){$A$}
\put(0.000000,6.353260){$K$}
\put(3.789854,6.353260){$G$}
\put(-0.450000,3.845684){$C$}
\put(-0.450000,2.953260){$D$}
\put(-0.450000,-0.2){$E$}
\put(1.4,3.2){$B$}
\put(2.0,3.02){$a$}
\put(1.789854,-0.350000){$F$}
\put(3.789854,-0.2){$L$}
\put(3.262354,-0.300000){$c$}
\put(3.889854,3.003260){$b$}
\endpspicture
}

\title{On the Notion of Equal Figures in Euclid}        
\author{Michael Beeson}
\address{ Michael Beeson\\
         San Jos\'e State University (emeritus) \\
        {\tt profbeeson@gmail.com}
       }        
\date{\today
}  

\begin{document}

\begin{abstract}Euclid uses an undefined notion of ``equal figures'', to 
which he applies the common notions about equals added to equals or 
subtracted from equals.  This notion does not occur in modern geometrical
theories such as those of Hilbert or Tarski.  Therefore to account
for Euclid in modern geometry, one must somehow replace Euclid's 
``equal figures'' with a defined notion.  In this paper we present
a new solution to this problem, and moreover we argue that ``Euclid 
could have done it''.   That is, it is based on mathematics that was
available in Euclid's time, including ideas related to Euclid's 
Proposition~I.44.  The proof uses the theory of proportions.  Hence we also discuss the ``early theory of proportions'', which has a long history.     
\end{abstract}
\maketitle

\section{Introduction}

The word {\em area} almost never occurs in Euclid's {\em Elements}, despite the 
fact that area is clearly a fundamental notion in geometry.%
\footnote{It does occur in English translation in Prop.~I.35, but in the 
context {\em parallelogrammic areas}, which according to Heath's commentary,
is intended mainly to emphasize that only four-sided figures are meant,
i.e., regular polygons of more than four sides with opposite sides parallel
are not meant.} Instead, Euclid speaks of ``equal figures.'' 
 Apparently a ``figure'' is 
a simply connected polygon, or perhaps its interior.  The notion is 
neither defined nor illustrated by a series of examples; for example, it 
is never made clear whether a figure has to be convex,  or even whether 
a circle is a figure, or whether a figure has an interior, or is just
made of lines.  

The notion of ``equal figures'' plays a central role in Euclid.  For 
example, the culmination of Book~I is the Pythagorean theorem.    
Nowadays we would, if required to express the theorem without algebraic
formulas, say that given a right triangle,
 the area of the square on the 
hypotenuse is the sum of the areas of the squares on the sides.  But 
Euclid said instead, that the square on the hypotenuse is equal to the
squares on the sides, taken together.  His proof shows how the two squares
can be cut up into  pieces that can be rearranged to make this equality of figures evident,
given earlier propositions about equal figures.  

Euclid did not define ``figure'',  and neither shall we.  For purposes of 
Euclid Book~I,  we can think of triangles and convex quadrilaterals, as these
are the only figures mentioned, but presumably Euclid did not mean {\em only} these,
but meant to include combinations of these as well.  Euclid also did not 
define ``equal figures'',  but simply treated it as a primitive (undefined) notion.
The main point of this paper is that he {\em could} have given a precise 
definition of ``equal figures'',  at least for equal triangles and convex quadrilaterals, and in that way proved all the theorems of Book~I without 
needing an undefined notion of ``equal figures.''

That Euclid did have area in mind when speaking of equal figures seems
clear from Book~II,  in which the whole thrust of the book is towards 
showing how, given any rectilineal figure, 
to find a square equal to the given figure;
one might interpret that as giving a method to calculate the area of any
rectilinear figure. 

\begin{figure}[ht]
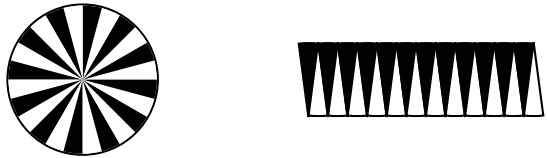

\caption{Archimedes proved the circle is equal to a certain rectangle, but 
he didn't use the word ``area.''}
\label{figure:archimedes}
\FigureArchimedesOne
\FigureArchimedesTwo
\end{figure}

Nor was Euclid alone in avoiding the word ``area.''
A century later, when Archimedes calculated the area of a circle,
he did not express his result by saying that the area of the circle is 
$\pi$ times the square of the radius.  Instead, he said that circle is 
equal to the rectangle whose sides are the radius and half the circumference.
(See Fig.~\ref{figure:archimedes}).  (So a circle did count as a figure for Archimedes!) %
\footnote{But Heron, in his {\em Metrica} of 50 CE \cite{heron-metrica},  did use the word 
\selectlanguage{greek}  ἐμβαδόν, \selectlanguage{english}
which is translated as ``area'',  and explains it this way:
``A cubit area is called when a square plot has each side of one cubit.''
The tile of {\em Metrica} shows Heron's concern with techniques for 
actually calculating areas; he even gives a numerical procedure for computing
approximate square roots.   He also has no problem multiplying four lengths
and then taking the square root to get an area.    But that answer was a number,
not a geometric length.  The identification of line segments with numbers 
was not a part of Greek mathematics.
}

Why did Euclid avoid the word {area}?  Not because he did not know 
that area can be measured;  it must have been for more abstract, mathematical 
reasons.   Let us consider his problem:
if he were to use the word, he would either have to {\em define} it, 
or put down some {\em postulates} about it.  Both choices offer some 
difficulties.  Area involves assigning a {\em number} to each figure, to 
measure its area.  It is therefore not a purely geometric concept.  
Moreover, even if one is willing to introduce numbers,  that just pushes
the problem back one step:  one must then define or axiomatize numbers.
Euclid knew that he did not know how to define area in general, so that
choice was out.  The other choice was to write down some axioms that area
obeys.  The most obvious one is additivity: if a figure can be cut into 
two pieces, then its area is the sum of the areas of the pieces.  But then 
there are delicate questions about the meaning of ``cut'' and ``piece.''
Euclid (or one of his unknown predecessors) 
discovered that it would possible to avoid all these complications
by replacing ``area'' by the concept of {\em equal figures}.  He noticed that if 
he used the word ``equal'',  and also re-interpreted ``taken together'' and ``taken from''
as if these operations were applicable to figures as well as to lines and angles,
then the additivity properties would look like special cases of the 
common notions  2 and 3, namely ``if equals be added to equals, the 
wholes are equal'', and ``if equals be subtracted from equals, the 
remainders are equal.''%
\footnote{This could now be explained using set theory, 
as   set-theoretic union and difference, but that is a development
only of the past century, and the verifications of Euclid's common notions for 
this notion still involve real numbers as well as sets.
}
  Common notion 5, ``the whole is greater than 
the part'',  could be taken to imply that a figure cannot be equal to 
a part of itself, and common notion 4, ``things which coincide with 
one another are equal to one another'',  could be interpreted to imply
that congruent figures are equal.  
Using these interpretations of the 
common notions, Euclid (thought he) could avoid all the complications mentioned above.%
\footnote{Actually,  Euclid needed one more property: halves of equal
figures are equal, used in Prop.~I.39.  The step that (implicitly) uses that property
occurs in Euclid's text without justification.}  

\begin{figure}[ht]
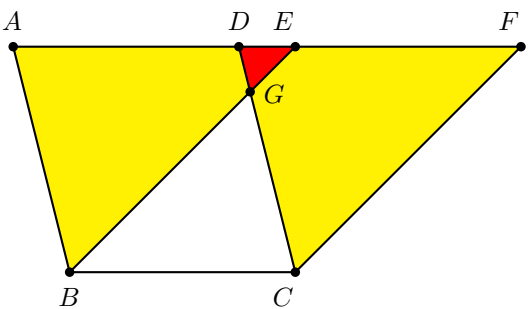

\caption{Euclid's proof of I.35}
\FigureOneThirtyFiveColored
\label{figure:I.35colored}
\end{figure}

We will give an example of how Euclid reasoned about equal figures,
namely  Euclid I.35.  See Fig.~\ref{figure:I.35colored}.

Euclid wants to prove the parallelograms $ABCD$ and $BCFE$
are equal.  He proves the triangles $ABE$ and $DCF$ are 
congruent.  Implicitly, he assumes $DEG$ and $DGE$ are
equal figures (that is, the order of listing the vertices 
does not matter).  Then ``subtracting equals from equals'', the 
yellow quadrilaterals are equal.  Then, ``adding equals to equals'',
he adds triangle $BCG$ (implicitly assuming $BCG$ is equal to $BGC$)
to arrive at the desired conclusion.  To formalize this 
proof, we needed so-called cut-and-paste axioms, as well as 
the axiom that $ABC$ and $ACB$ are equal triangles.  In this paper,
we will show how to define ``equal figures'' so these propositions
can be proved instead of assumed.

Although it appears to the modern eye (e.g. \cite{hartshorne}) that 
Euclid meant ``figures with equal area'' when he said ``equal figures'',
it is worth noting that not only does he never mention the word ``area'',
but he also never speaks of one figure being greater than another, 
although certainly areas can be compared.  He never applied the common 
notions that mention ``greater than'' to figures.  Perhaps he thought
``$A$ greater than $B$'' generally means that $B$ is equal to a part of $A$;
that definition, if applied to figures, does not lead to the same 
laws that ``greater than'' enjoys for lines.  This line of thought 
casts a shadow of doubt on the theory that ``equal figures''
meant ``equal area'', without suggesting another interpretation.  

Later generations of mathematicians were not willing to accept Euclid's
over-liberal interpretation of the common notions in support of 
``equal figures.''  See the summary 
discussion with many references on pp.~327--328 of \cite{euclid1956}.  In
particular, once mathematicians had some experience with axiomatization,
it became obvious that ``equal figures'' is not a special case of 
equality, since equal figures cannot be substituted for each other 
in every property.   Instead, it is a new relation, and the original
choice that Euclid finessed faces us directly: we must either define
or axiomatize the notion.

Euclid does not mention ``equal figures''   
until Prop.~I.35.  The reader is urged to look at the proofs
of Euclid's Propositions I.35, 42, 43, 47, and 48,  to identify
the lines of Euclid's proofs where common notions about equal figures
are used.    They are relatively few in number, but crucial 
to Euclid's development.  For example, Euclid Book~I culminates in 
the Pythagorean theorem, which cannot even be stated without the 
notion of equal figures.   

We mention a few matters of notation.  Euclid usually used upper-case
letters for points.  We use both upper-case and lower-case letters for points;
this allows the notation to suggest correspondences, as in ``Triangle $ABC$ 
is congruent to triangle $abc$.''  Since the force of tradition in geometry
is strong, we justify this choice:  One might also choose to use primes or 
stars or circumflexes, as in ``$ABC$ is congruent to $A^\prime B^\prime C^\prime$.''
Such notations are hard to fit into diagrams, and cannot be cut-and-pasted into
computer systems for formal proofs, leading to errors of transcription.  

A second matter of notation is ``betweenness''.  Starting with Pasch \cite{pasch1882},
geometers have used the relation ``$B$ is between $A$ and $C$''; 
Hilbert \cite{hilbert1899} introduced the notation $\B(A,B,C)$, which we shall
use in a few places.  Following Hilbert we take it to mean ``strict betweenness'';
that is, it means $B$ lies on the line (segment) $AC$ and is not equal to $A$ or $C$.
\medskip

\noindent
{\em Acknowledgements}. I would like to thank Victor Pambuccian and 
Vincenzo De Risi for their help with the history of this subject.  I would
like to thank the anonymous referee for an extremely careful reading and 
many helpful suggestions.

\section{Possible ways to define area and equal figures}

The problem we are trying to solve is to eliminate the 
need for a primitive notion of ``equal figures.'' 
  In this section we explain the 
principal approaches to this problem that have been tried.

\subsection{Hilbert's equidecomposition}
 When Hilbert wrote his influential book
\cite{hilbert1899}, he chose to define the notion.  His definition 
was still used in the much more modern book \cite{hartshorne}, p.~200.
The problem with that definition is that it mentions the notion of 
natural number,  in speaking of cutting a figure into ``a finite number''
of triangles.  It is therefore not a purely geometric notion.  More
technically, it cannot be expressed in a first-order language with 
only geometric variables.  This concept does not permit us to 
define {\em area}, but only {\em equal area} for polygons. 
That would, in itself, not be a problem, as Book~I needs only triangles and 
convex quadrilaterals, but Hilbert himself showed that (even for triangles)
the claim that equal figures are equidecomposable requires Archimedes's axiom
\cite[\S 19]{hilbert1899}.

\subsection{Defining area by calculus}
Since the eighteenth century, we have had
 the option to define area using integrals,
an option that was not available to Euclid, who wrote even before 
Archimedes's work on the circle pioneered the use of limits.  But that 
too is not a purely geometric definition,  since it involves real numbers
and functions as well as limits.  However, it is the {\em only}
workable definition of area that mathematicians have found.

\subsection{Descartes's geometric arithmetic}
 Descartes and later Hilbert 
defined {\em geometric arithmetic}.  That is,  we fix
a certain line (the $x$-axis) and two points on that line (0 and 1).
Then certain geometric 
constructions exist for defining the addition and multiplication of 
line segments with $0$ as one endpoint,  and one can give geometric 
proofs of the laws of arithmetic, such as the associative, commutative,
and distributive laws. 

When Hilbert, and later Tarski, worked on the formalization of elementary
geometry, their principal aim was to develop this theory of ``geometric 
arithmetic.''  They did not follow Euclid and did not, for example, prove
Prop.~I.35 and subsequent propositions that depend on equal figures, but 
proceeded by the most direct route to geometric arithmetic. That route led
through the theorems of Pappus and Desargues.  

It is not the case that just because we can multiply we can define
areas, even of triangles.  The problem of interpreting what Euclid
meant by ``equal area'' is not automatically solved by defining
geometric arithmetic.    

Besides, segment arithmetic was conceptually alien to Greek thought. 
The Greeks never multiplied
lengths to get lengths.   A length times a length produced a rectangle;
multiplying three lengths produced a solid;  four lengths could not be 
multiplied.  Even Heron, who had no problem multiplying numbers to 
compute areas, did not think he was multiplying lengths to get lengths.
That conceptual 
step was not taken until Descartes, in 1637.   Even Vieta,  who
introduced  using letters for quantities in algebra,
still adhered to the ``principle of homogeneity'' in which all algebraic 
terms in an equation had to have the same degree.%
\footnote{  For example, see \cite{vieta1646}, p.~86,  where Vieta writes
A cubus + B quad. in A, equetur B quad.in Z,  or in modern symbols,
$A^3 + B^2A = B^2Z$ instead of $x^3 + px = q$.  See \cite{hartshorne-vieta}
for further discussion.  Incidentally, one sees both Vi\`ete and Vieta, 
the French and Latin spellings of the name.}  

\subsection{Axiomatize the properties of area}
We could try to axiomatize the properties of an ``area function'' 
from sets of points to numbers.  To avoid the complications of ``numbers'',
Hartshorne tried this approach using any ordered Abelian group for the
values of area.  See pp.~205ff. of \cite{hartshorne} , where it is shown how 
to do this using a minimum of assumptions.   But this is also 
not a geometric notion, as the area function has to take values 
{\em somewhere}.  Hartshorne uses any abelian group.

\subsection{Direct axiomatization of the equal-figure axioms}

The definitions of ``equal figure'' discussed above all 
are unsatisfactory, since they require the concepts of real number,
or natural number, or both; or else, in the case of segment arithmetic,
do not actually lead to a definition of equal area.

When we set out, in previous work \cite{beeson2019}, to
formalize Euclid Book~I,  we saw no 
other alternative.   Following Sherlock Holmes's maxim that when the impossible is
eliminated, what remains is the truth,  we chose the remaining
alternative:  to axiomatize the notion.   
  That was also implicitly the choice 
in \cite{hartshorne}, where Hilbert's definition is only used for long
enough to establish the properties that Euclid used, and to it is 
added ``de Zolt's axiom'', that if $Q$ is a figure contained in another
figure $P$, and $P-Q$ has a nonempty interior, than $P$ and $Q$ are not
equal figures.  Hartshorne says that he does not know a purely geometrical
proof of this from Hilbert's definition of equal figures; probably he 
means that to prove it, you must prove that two figures are equal in 
Hilbert's sense if and only if they have equal area, in the sense of 
area defined by integrals. 
 
 The common notions of Euclid include 

\begin{itemize}
\item CN2:
If equals are added to equals, then the wholes are equal.
\item CN3:
If equals are subtracted from equals, then the remainders are equal.
\item CN5:
The whole is greater than the part.
\end{itemize}

Euclid meant these to apply to various sorts of ``things''; in Book~I
he used them for triangles and quadrilaterals, which are the only ``figures''
in Book~I.   
Our plan in \cite{beeson2019}, and Hartshorne's plan in \cite{hartshorne},
was to translate the ``common notions'' mentioned above into first-order
axioms.  Since the variables range over points,  triangles are just 
triples of points,  so ``equal triangles'' is a 6-ary relation, 
$ET(A,B,C,a,b,c)$.  Then we need an 8-ary relation $EF$ for 
``equal quadrilaterals''.  ($EQ$ was already in use for equality, 
so we used $EF$ for ``equal figures.'')  In Euclid Book~I, only 
triangles and quadrilaterals are used, so we stopped there. 
The parts of a figure are smaller figures.  New figures are
constructed, in Euclid's proofs, by cutting off triangles from larger
triangles or quadrilaterals,  and also by pasting on figures to other 
figures along a common edge. 

Euclid also used two principles about equal figures without ever
formulating them as axioms or common notions:  halves of equals are equal,
and doubles of equals are equal.  These also correspond to axioms in 
our list of equal-figures axioms. 

 Since there are several ways to do this 
cutting and pasting, we get several ``equal-figures axioms'' this way.
The last common notion, CN5,  becomes ``de Zolt's axiom'', after 
the person who first formulated it. 
These axioms express the ways in which Euclid used the ``common notions''
as applied to figures. 

This approach to the notion of equal-figures by direct axiomatization
was successful:  we constructed proofs ``faithful to Euclid'' 
that were verifiably logically correct.   
Nevertheless,  the addition of fifteen new axioms to Euclid is 
a bit unsatisfying,  even if they do correspond well to Euclid's actual 
proofs.  What we really want is a treatment of geometry that 
\smallskip

(i)  Corresponds as well as possible to Euclid (but corrects the errors), and 
\smallskip

(ii) Has axioms that correspond well to basic geometric intuitions, and 
\smallskip

(iii) meets modern standards of rigor.
\smallskip

The system of axioms that we used in \cite{beeson2019},
or the similar system in \cite{hartshorne},
 Euclid meets these standards, if 
you think that the fifteen equal-figure axioms correspond to a basic 
geometric intuition.  But it seems that they do not really:  they are justified
by an appeal to our intuitions about area,  and area cannot be expressed
purely geometrically, as it fundamentally involves using numbers to measure area.

\section{Aim and methods of this paper}

The contribution of this paper is to eliminate the ``equal figures'' axioms
by {\em defining} the notions of ``equal triangles'' and 
``equal quadrilaterals'', by a definition that Euclid could have given,
and proving the properties expressed in the ``equal figures'' axioms,
so that Euclid Books~I to IV could be developed without the 
equal figures axioms.

This aim would not be met by following Hilbert and Tarski, 
who first define segment arithmetic and then use 
it to define area.  That route would  take us 
 far from Euclid, who never thought of multiplying two 
line segments.    While it would
meet the requirement to eliminate the ``equal figures'' axioms by defining that
notion, it would not meet the requirement that ``Euclid could have done it.''
 
Our method to achieve this aim
 is to define ``equal rectangles'' using a figure much like the one 
Euclid uses for Prop.~I.44.
and use that to define ``equal triangles'' and ``equal quadrilaterals''  and
use those defined notions to prove the propositions of Euclid Book~I.
To reiterate:  we will
\begin{itemize}
\item  define ``equal rectangles'' and ``equal triangles'', and
\item prove the propositions of Euclid Book~I using those defined
notions (rather than extra axioms), and 
\item  use proofs ``in the spirit of Euclid''
\end{itemize}
\medskip

It follows that the equal-figures axioms are 
actually superfluous, in the sense that, using the new
definition of ``equal figures'',  we could formalize Euclid Book~I
directly, without adding any equal-figures axioms.  But we 
then take one step more, and show that the equal-figures axioms
can in fact all be proved. 

Our proof that the defined notion
of ``equal figures'' has the required properties 
uses  some theorems about similar triangles and proportion.
  In Euclid, those theorems are present, but only in Book V,  after
the development in Book IV of Eudoxes's theory of ``magnitudes.'' 
To carry out our program thus requires the demonstration that 
Euclid could have developed the necessary theory of proportion 
without using the Axiom of Archimedes (on which Book V depends),
and preferably without using Book III (theorems about circles),
and of course without using ``equal figures,''
so that it would be available in the last third of Book~I.   
 It turns out that we are not the first to seek an earlier development of 
the theory of proportion; in \S\ref{section:kupffer2} below, we 
discuss this subject at length, with historical notes.

The proofs that we give in this paper are {\em informal} in the 
sense that they have been written for humans to read,  but 
{\em rigorous} in the sense that they can be carried out in 
any formal system adequate for elementary geometry.  For example
in Hilbert's system, or Tarski's system, or the system used in 
\cite{beeson2019},  or the textbook \cite{hartshorne}, or
(apart from the errors corrected in \cite{beeson2019}) in Euclid's own
system, minus steps about equal figures that we are trying to justify here.%
\footnote{It is therefore {\em not} necessary to fix a particular first-order
version of Euclid to check this paper, unless of course, one wants to check
the proofs by computer.  In that case, refer to \cite{beeson2019} or the 
perhaps more accessible \cite{beeson2022}.}
This paper could have been written and read in the nineteenth century,
before the invention of modern logic and computers.  It would, of course,
have needed to be after the theorems on proportionality in \S~\ref{section:kupffer}.

\section{Similar triangles and proportion}
\label{section:proportionality}

One of the principles about ``equal figures'' that Euclid uses
in at least two crucial places is that $ABC$ is equal to $BCA$. 
Under the definition of ``equal triangles'' that we give below, 
this principle turns out to require the basic theorems about 
similar triangles and proportion for its proof.  Specifically,
if two right triangles have their hypotenuses and one leg proportional,
then they have corresponding angles equal.  
 Also, the proof that ``equal rectangles''
is a transitive notion seems to require properties of proportion.

In this section, we state the definitions and properties of 
proportion that we will use.  Euclid proved these theorems in
Book VI, using the results of Book~V, which is based on the 
axiom of Archimedes. Because of the reliance on Archimedes's axiom
and possibly the use of equal figures, this is not useful for us.
But Paul Bernays proved these results without using Archimedes's
axiom, and by means acceptable for our purposes,  in his 
Supplement~II to \cite{hilbert1899}.  
In this section, we will give the definition of proportion,   prove
some easy lemmas,  and make it clear that there are exactly 
two non-trivial theorems in the subject, namely the ``interchange theorem''
and the ``fundamental theorem'',  which are stated below, and 
proved by Bernays and Kupffer  \cite{kupffer1893, schur1902}.  
Their proofs and some relevant history will be discussed in \S\ref{section:kupffer}.

\begin{figure}[ht]
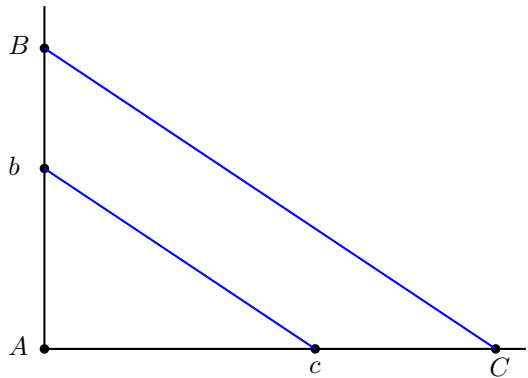

\caption{Here $AB:AC = Ab:Ac$ because $BC\  ||\  bc$.}
\label{figure:proportiondefn}
\FigureProportionDefn
\end{figure}

\begin{definition}\label{definition:proportion}
$AB:AC = Ab:Ac$ if  $BAC$ is a right triangle, $b$ is 
on ray $AB$, $c$ is on ray $AC$, and $BC$ is parallel to $bc$
(or $BC = bc$),
as shown in Fig.~\ref{figure:proportiondefn}. 
More generally $PQ:RS = pq:rs$ if $PQ$, $RS$, $pq$, and $rs$ are
congruent to such segments $AB$, $AC$, $Ab$, and $Ac$.
\end{definition}

Formally, we have defined a relation taking eight arguments
of type ``point''.%
\footnote{We have not defined $AB:AC$ as a function taking four 
points, or two segments; the use of the equality symbol and colon in 
informal writing is just an abbreviation for the 8-argument 
relation.  Bernays in his Supplement~II to \cite{hilbert1899}
does define $a:b$ for segments $a,b$ to be, in effect, the 
angle whose tangent is $b/a$; but he never makes any use of that
definition other than to verify that the equality of such 
angles implies the definition of proportionality 
we give here.} 

It follows from the symmetry of ``parallel'' and congruence that  $pq:rs = PQ:RS$ 
if and only if $PQ:RS = pq:rs$, which we will use without further explicit mention.  We also have the following simple property:

\begin{lemma} \label{lemma:flip} 
$pq:rs = PQ:RS$ if and only if $rs:pq = RS:PQ$.
\end{lemma}

\noindent{\em Proof}.  Immediate from the definition of proportion 
and the symmetry of ``parallel''.
\smallskip

The following definition is Euclid's wording, but the meaning
he attached to ``proportional'' is different.

\begin{definition}\label{definition:similar} 
Two triangles $ABC$ and $abc$ are similar
if their corresponding angles are equal and their corresponding sides 
are proportional, that is, $AB:AC = ab:ac$.  
\end{definition}

\noindent{\bf Notation}.  If we write $a:b = c:d$, the letters
$a$, $b$, $c$, $d$ stand for pairs of point variables.  We may think of 
these pairs as segments, but formally the relation $a:b = c:d$  is an 8-ary 
relation on points.  That the abbreviated notation is convenient is 
illustrated by shortening the statement of Lemma~\ref{lemma:flip}: 
$u:v = p:q$ if and only if $v:u = q:p$.  Variables occurring next to 
a colon cannot be points; they must be pairs of points,  so there is no 
ambiguity in using lower-case letters in this way, as well as for points.

\begin{lemma}[transitivity]\label{lemma:proportion-transitive}
If $a:b = c:d$ and $c:d = e:f$ then $a:b = e:f$.  
\end{lemma}

\noindent{\em Proof}.  Immediate from the transitivity of ``parallel.''
\medskip

Bernays proves the existence and uniqueness of the fourth proportional;
we present versions of his proofs, 
 to check that they work in the present framework, that is, 
without relying on parts of Euclid past I.35. 

\begin{lemma}[existence of the fourth proportional] 
\label{lemma:proportion4e} 
For each $p,q,r$, there exists a segment $x$ such that $p:q = r:x$
\end{lemma}

\noindent{\em Proof}.  Given segments $p,q,r$, we let $A$ and $B$
be the endpoints
of $p$.  Erect a perpendicular $AC$ to $AB$ at $A$, with $AC = q$.
Then construct point $b$ on ray $AB$ so that $r = Ab$.  
 We may assume $b \neq B$, since if $b=B$ then $p=q$, so we may take $x=r$.
Then we may construct
line $\ell$ through $b$ parallel to $BC$.  (See Fig.~\ref{figure:proportiondefn}; we have to prove point $c$ in the 
figure exists.) Since $BAC$ is a right
angle, $A$ does not lie on $BC$, and angle $ACB$ is less than a right angle.
Then by Euclid 5,  line $AC$ (possibly extended) meets line $\ell$ 
in a point $c$.  Then by the definition of proportionality,
we have $AB:AC = Ab:Ac$; taking $x = Ac$ we have $p:q=r:x$.
That completes the proof. 

\begin{lemma}[uniqueness of fourth proportional] \label{lemma:proportion4}
For each $a,b,c$, there is exactly one segment $x$ such that $a:b = c:x$
(up to congruence).
\end{lemma}

\begin{figure}[ht]
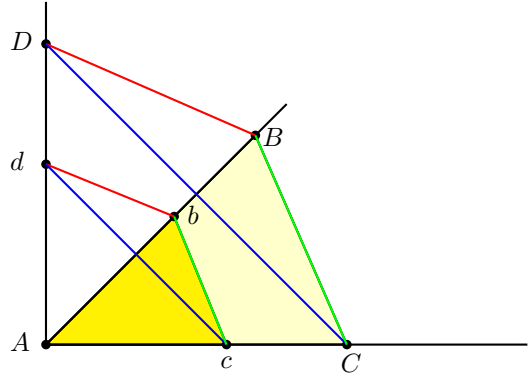

\caption{If $BC \ || \ bc$  
 then $AB:Ab = AC:Ac$}. 
\label{figure:proportion1}
\FigureProportionOne
\end{figure}

\noindent{\em Proof}.  This is an easy consequence of the definition.  
We spell out the details to facilitate formalization.   Suppose $a:b = c:d$ and $a:b = c:x$.   By 
Lemma~\ref{lemma:proportion-transitive}, $c:d = c:x$.    Since proportion is defined up 
to congruence, we may assume that in Fig.~\ref{figure:proportion1},  $AB$  and $Ab$ are both 
equal to $c$, while $AC$ is $d$ and $Ac$ is $x$.  But then,  $b=B$, so by the definition
of proportion, $BC$ and $bc$ coincide.  Hence $c=C$, as both are the intersection point
of line $BC$ with line $AC$.  Hence $Ac=AC$; but that is $x=d$.  That completes the proof.

The following theorem is Euclid V.16, except that Euclid's
definition of proportion is not the same.

\begin{theorem}[Interchange theorem]
\label{theorem:proportion-interchange}
 If $a:b = p:q$, then $a:p = b:q$.
\end{theorem}

This theorem will be proved in \S\ref{section:kupffer}.

\begin{lemma} \label{lemma:proportional-legs} If two right triangles
have their legs proportional, then their corresponding angles are equal.
\end{lemma}

\noindent{\em Proof}. Since proportionality is defined up to congruence,
we may assume without loss of generality that the two right 
triangles $ABC$ and $aBc$ have their right angle at $B$ in common,
and $a$ lies on ray $BA$ and $c$ lies on ray $BC$.  Then by the definition
of proportionality, $AC$ is parallel to $ac$ (or coincident).  Then 
the corresponding angles are either identical, or are corresponding 
angles of the traversals $BA$ and $BC$ or the parallel lines $AC$ and $ac$.
Hence the corresponding angles are equal.  That completes the proof 
of the lemma.

The following theorem was named by Bernays in
his Supplement II to \cite{hilbert1899}.

\begin{theorem}[fundamental theorem of proportion] 
\label{theorem:fundamental}
If two parallels delineate the segments $AC$ and $Ac$ on one 
side of an angle and $AB$, $Ab$ on the other side of the angle,
then $AB:Ab = AC:Ac$.
\end{theorem}

This theorem will be proved in \S\ref{section:kupffer}.
Fig.~\ref{figure:proportion1} illustrates the fundamental theorem.
In the figure, $AD = AC$ and $Ad = Ac$. 
Modulo the interchange theorem, we could write the conclusion
$AB:AC = Ab:Ac$, which is equivalent to $DC \ || \ dc$ in the
figure, by the definition of proportionality, where

The figure makes it clear that the fundamental theorem is a 
consequence of Desargues's theorem.  Whether it implies Desargues's 
theorem in some simple way we do not know. Bernays 
proved it by much more elementary means, discussed in \S\ref{section:kupffer}.

\begin{corollary} \label{lemma:proportion1}  If two triangles $ABC$ and $abc$
have corresponding angles equal, then $AB:AC = ab:ac$.
\end{corollary}

\noindent{\em Proof}. Since proportionality is defined up to 
congruence, we may assume $a=A$ and $AB$ and $Ab$ are collinear
and $AC$ and $Ac$ are collinear. Then the fundamental theorem applies.
That completes the proof.

\begin{corollary}\label{lemma:proportional-implies-similar} If two right triangles have their hypotenuses and one leg proportional, then their
corresponding angles are equal.   
\end{corollary}

\begin{figure}[ht]
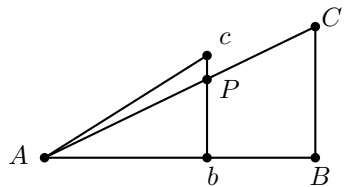

\caption{$AB:ab = AC:ac$ implies angles $CAB$ and $cAb$ are equal.}
\FigureProportionalSimilar
\label{figure:proportionalsimilar}
\end{figure}

\noindent{\em Proof}. Please refer to Fig.~\ref{figure:proportionalsimilar}.
Since proportionality is defined up to congruence,
we may assume without loss of generality that the two right 
triangles are $ABC$ and $Abc$, with vertex $A$ in common and
$b$ on ray $AB$, and right angles at $b$ and $B$, 
and $AB:ab = AC:ac$.  We must show that
$c$ lies on ray $AC$.  

We have two straight lines $Ac$ and $Bc$,
and line $AB$ falling on those two straight lines.  The interior angles 
$cAb$ and $Abc$ are together
less than two right angles, since $Abc$ is right and $cAb$ is less than right, since $Abc$ is right.
Therefore, by Euclid 5, $AC$ meets $bc$ in a point $P$ on the same side of $AB$ as $c$.  
Then we have
\begin{eqnarray*}
AB:ab &=& AC:ac \mbox{\qquad\ \ by hypothesis} \\
AB:ab &=& AC:AP \mbox{\qquad by the fundamental theorem, since $BC \ || \ bc$}\\
Ac &=& AP    \mbox{\qquad\qquad\ by the uniqueness of the fourth proportional} \\
\end{eqnarray*}
It remains to show $c=P$, which seems visually obvious, but
is not trivial to prove.
Assume, for proof by contradiction, that $\B(b,P,c)$. 
Then triangle $AcP$ is isosceles, so angle $AcP$ is equal 
to angle $aPc$.  Since the angles of triangle $AcP$ make together
two right angles, angle $APc$ is less than a right angle.  But
also angle $APb$ is less than a right angle, since angle $AbP$ 
is a right angle.  But angles $APb$ and $APc$ are supplements,
since $\B(b,P,c)$, so together they make two right angles, contradiction.
Hence $P$ is not between $b$ and $c$.  

Now assume, for proof by contradiction, that $\B(b,c,P)$. Again
triangle $AcP$ is isosceles, so angle $AcP$ is less than a right angle.
Also angle $Abc$ is less than a right angle.  But now 
angles $Abc$ and $AcP$ are supplements, contradiction.

Hence neither $\B(b,c,P)$ nor $\B(b,P,c)$.  But $P$ lies on line $bc$; 
therefore $P=c$.  Then $c$ lies on ray $AC$.  Hence angle 
$cAb$ is equal to angle $CAB$.  Hence also angle $Acb$ equals angle
$ACB$.  That completes the proof of the lemma. 
\medskip

We next give two ``combination rules'' introduced and proved
by Bernays, \cite[p.~204]{hilbert1899}. 
 We repeat the proofs
given by Bernays, as they are short.  The reader can check that they 
use only the interchange theorem and the existence of the fourth proportional.

\begin{lemma}[Bernays's ``first combination rule''] \label{lemma:combination1}
\ \\ If $a:b=a^\prime:b^\prime$  and $b:c = b^\prime:c^\prime$ then
$a:c = a^\prime:c^\prime$.
\end{lemma}

\noindent{\em Proof}.   Suppose 
If $a:b=a^\prime:b^\prime$ and $b:c = b^\prime:c^\prime$.
Then 
\begin{eqnarray*}
a:a^\prime = b:b^\prime  &&\mbox{\qquad by the interchange theorem}\\
b:b^\prime = c:c^\prime &&\mbox{\qquad by the interchange theorem}\\
a:a^\prime = c:c^\prime && \mbox{\qquad by Lemma~\ref{lemma:proportion-transitive}}\\
a:c = a^\prime:c^\prime &&\mbox{\qquad by the interchange theorem}
\end{eqnarray*}
That completes the proof.

\begin{lemma}[Bernays's ``second combination rule''] \label{lemma:combination2}
\ \\ If $a:b=b^\prime:a^\prime$ and $b:c = c^\prime:b^\prime$ then
$a:c = c^\prime:a^\prime$.
\end{lemma}

\noindent{\em Proof}. By Lemma~\ref{lemma:proportion4e}, 
let $u$ be the fourth proportional to $a,b,c^\prime$. 
\begin{eqnarray*}
a:b = c^\prime:u   &&\mbox{\qquad since $u$ is the fourth proportional}\\
a:b = b^\prime:a^\prime  &&\mbox{\qquad by hypothesis} \\
c^\prime:u  = b^\prime:a^\prime && \mbox{\qquad by the preceding two lines}\\
c^\prime:b^\prime = u:a^\prime && \mbox{\qquad by the interchange theorem}\\
b:c = c^\prime:b^\prime &&\mbox{\qquad by hypothesis} \\
u:a^\prime = b:c  && \mbox{\qquad by Lemma~\ref{lemma:proportion-transitive}}\\
c^\prime:u = a:b &&\mbox{\qquad since $a:b=b^\prime:a^\prime = c^\prime:u$}\\
c^\prime:a^\prime = a:c && \mbox{\qquad by Lemma~\ref{lemma:combination1}}
\end{eqnarray*}
That completes the proof.

We have now derived the theory of proportionality from just two 
theorems, the interchange theorem and the ``fundamental theorem
or proportionality.  
In \S\ref{section:kupffer}, we will discuss Bernays's  and Kupffer's
proofs
of those theorems, with due
attention to the extensive history of this subject.  In the 
meantime, we move on to develop our theory of equal figures, 
using these facts about proportion.

\section{Defining equal rectangles and equal triangles}
\subsection{Order matters} \label{section:orientation}
Is ``triangle $ABC$'' the same triangle as ``triangle $BCA$''? 
If you are willing to accept the answer ``no'',  you may skip this subsection;
otherwise read on. 

Neither Euclid nor his modern successors defined ``triangle''; 
triangles are always introduced and referred to by triples of points $ABC$.
In Proposition~I.35,  where equal figures are used,  it is taken for 
granted that $DGE$ is equal to $EGD$, so that ``subtracting'' these two 
triangles from larger ones is ``subtracting equals.''  But Euclid 
never commits himself to saying that $DGE$ and $EGD$ are the {\em same}
triangle, or are {\em different but equal} triangles.   

We can't tell which he meant, since in Euclid,
 triangles are just given by {\em ordered}
triples of points.   Since we do not actually have triangles in our
ontology, perhaps it doesn't matter much.  But we do define
a relation of ``triangle congruence.''  This is a 6-ary relation
on points, defined in terms of segment congruence (more formally,
the ``equidistance relation'').  Namely,
 $ABC$ is congruent to $abc$ if $AB$ and $ab$ are congruent,
$BC$ and $bc$ are congruent, and $AC$ and $ac$ are congruent.  
According to this definition, ``order matters'':   $ABC$ will 
be congruent to $BAC$ only when $BC$ and $AC$ are congruent.

Actually, Euclid never defines ``congruent triangles.'' 
For example, look up his statement of I.4 (the SAS criterion). 
But our definition is  in accordance with a long tradition since 
Euclid:  when we say that $ABC$ and $abc$ 
are congruent triangles, we intend to imply that $AB$ and $ab$ are 
corresponding sides,  $BC$ and $bc$ are corresponding sides, $ABC$ and $abc$ are
corresponding angles, etc.   That is,  the order in which we mention
the vertices of triangle $ABC$ does matter, if we are to speak of congruent triangles
and have the names of the triangles convey which are the corresponding
sides and angles in the congruence.%
In other words:  this is not a deep philosophical issue about the nature of 
triangles.  It is a convention concerning how we describe triangles by mentioning
their vertices.  Whether $ABC$ is ``the same triangle as $BAC$'' or 
``a different triangle from $BAC$'' simply never comes up, because
we don't say what a triangle is. 

We define the {\em base of triangle $ABC$}  to be $AB$.  This is 
just a way of selecting the first two of the three points.  

\subsection{Equal rectangles}
We begin by defining the 8-ary relation {\tt ER}, 
 ``equal rectangles.''   Two rectangles are 
{\em congruent} if their sides are pairwise equal.  Explicitly,
$ABCD$ is congruent to $abcd$ if $AB = ab$ and $BC = bc$.

\begin{figure}[ht]
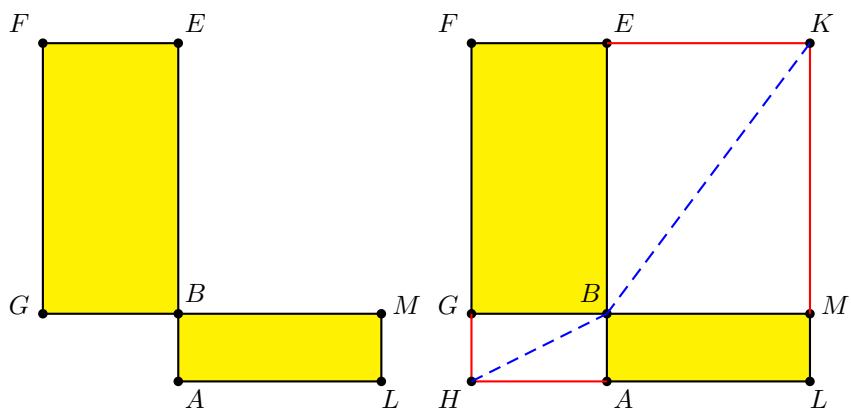

\center{\EqualRectanglesFigureOne \EqualRectanglesFigureTwo}
\caption{(Left) Place copies of two given rectangles as $BEFG$ and $\mbox{\it BMLA}$
with two sides collinear. (Right) 
The other sides (extended) meet by Euclid 5 forming a large rectangle.
Then $BEFG$ and $\mbox{\it BMLA}$ are defined to be equal if $\B(H,B,K)$.
In the case shown, they are not equal.  
}
\label{figure:equalrectangles}
\end{figure}

\begin{definition} \label{definition:equalrectangles}
Any two given rectangles $R$ and $S$ are congruent to rectangles placed like 
$FEBG$ and $BMLA$
in Fig.~\ref{figure:equalrectangles}, where $\B(A,B,E)$ and 
$\B(G,B,M)$, and $ABGH$ and $\mbox{\it BMKE}$ are rectangles.
  By definition the two 
rectangles $BEFG$ and  $BMLA$ are {\bf equal rectangles} if and only if $\B(H,B,K)$,  that is, 
the line $HK$ passes through the common vertex $B$ of the two rectangles.
We then say that  
the two original rectangles $R$ and $S$ are equal. 
\end{definition}

\begin{figure}[ht]
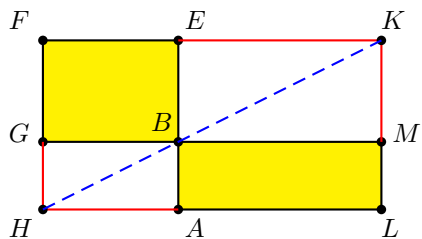

\center{\EqualRectanglesFigureThree}
\caption{In this case the two rectangles are equal,
since $B$ is between $H$ and $K$.  That is, the 
two dashed lines form one straight line.}
\label{figure:equalrectangles3}
\end{figure}
\FloatBarrier

The connection between equal rectangles and the theory of proportion
is given in the following lemma:

\begin{lemma} \label{lemma:ERproportion} The rectangle with base $b$ and height $a$
is equal to the rectangle with base $c$ and height $d$  if and only if 
$b:c = d:a$, and also if and only if $b:d = c:a$.
\end{lemma}

\noindent{\em Proof}.  The two proportionality statements are 
equivalent, by Theorem~\ref{theorem:proportion-interchange}. 
Suppose rectangle 1 has base $b$ and height $a$, and rectangle 2
has base $c$ and height $d$.  

Please refer to Fig.~\ref{figure:equalrectangles3}, in which 
rectangle 1 is  $FEBG$ and rectangle 2 is $BMLA$.
 
Suppose the two rectangles are equal.  Then $B$ is between $H$ and $K$,
so $HAB$ and $BMK$ have their corresponding angles equal.
By Corollary~\ref{lemma:proportion1}, $KM:BM = AB:AH$.
By the interchange theorem, $KM:AB=BM:AH$. Now $KM = EB = a$,
$AB = d$, $BM = AL = c$ and $AH = GB = b$.  Thus $a:d = c:b$.
By Lemma~\ref{lemma:flip} then $d:a = b:c$.  That completes the left-to-right
direction of the proof.  

Now suppose $d:a = b:c$.  The last few steps are reversible,
leading to $KM:BM=AB:AH$.  Then by Corollary~\ref{lemma:proportional-legs}, the right angles $HAB$ and $BMK$ have their corresponding
angles equal.  I say that implies $B$ is between $H$ and $K$.
Indeed, since $GM$ is parallel to $HL$, the extension of line $HB$
through $B$ makes an angle with $BM$ equal to angle $KBM$.  
If we knew that $K$ is on the opposite side of $GM$ from $H$,
we could conclude by Euclid~I.7 that 
$K$ lies on that extension; that is, $B$ is between $H$ and $K$.%
\footnote{Euclid~I.7 is quite a bit more difficult to prove 
than Euclid thought; it is hard to prove that an angle cannot be both 
equal to and less than another angle.  Hilbert avoided the difficulty
by including uniqueness in his angle-copying axiom.}

To prove that $K$ is on the opposite side of $GM$ from $H$,
we note first that $H$ is on the opposite side of $GM$ from $F$
since $\B(H,G,F)$.  So it suffices to prove that $H$ is on the 
same side of $GM$ as $K$.  For that it suffices to show that 
$FA$ and $AK$ both meet $GM$.  Those assertions are applications
of a lemma called {\tt parallelpasch} that we proved during 
our formalization of Euclid Book~I.  The lemma says that if 
a point $A$ lies on the extension  of one side $EB$ of a parallelogram
$EBGF$, then $AF$ meets $BG$; that is exactly what we need here.

That completes the proof of the lemma.

\begin{lemma} \label{lemma:ERequivalence}
``Equal rectangles'' is an equivalence relation.
\end{lemma}

\begin{figure}[ht]
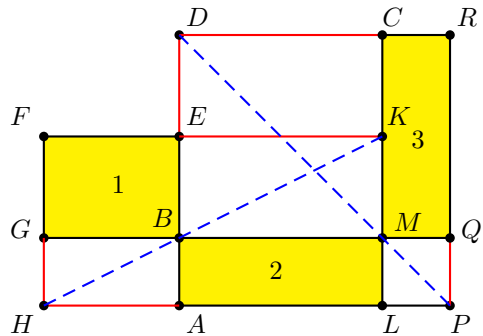

\caption{If rectangle 1 equals rectangle 2 equals rectangle 3, then 
rectangle 1 equals rectangle 3.}
\label{figure:ERtransitive}
\FigureERTransitive
\end{figure}

\noindent{\em Proof}.  We omit the proofs of reflexivity and symmetry.
The difficult part is transitivity.  Suppose rectangle 1 is equal 
to rectangle 2,  and rectangle 2 is equal to rectangle 3.  We 
have to prove that rectangle 1 is equal to rectangle 3.   See 
Fig.~\ref{figure:ERtransitive}.   We have
\begin{eqnarray*}
MQ:BA = BM:MC &&\mbox{\qquad since $3 = 2$} \\
BA:BG = BE:BM   &&\mbox{\qquad since $2 = 1$} 
\end{eqnarray*}
We apply Bernays's second combination rule (Lemma~\ref{lemma:combination2});
with
\begin{eqnarray*}
a = MQ, b = BA, b^\prime = BM, a^\prime = MC,
c = BG, c^\prime = BE
\end{eqnarray*}
The hypotheses of the combination rule become the above equations.
The conclusion is $a:c = c^\prime:a^\prime$; that is,
$$ MQ:BG = BE:MC.$$
Then rectangle 3 is equal to rectangle 1, by Lemma~\ref{lemma:ERproportion}.
That completes the proof of the lemma.

\begin{lemma} \label{lemma:ER1} Rectangles $ABCD$ and $BCDA$ are
equal; also for any permutation of the vertices that is still a rectangle.
\end{lemma}

\begin{figure}[ht]
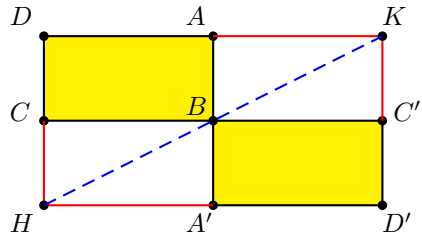

\caption{$ABCD$ is equal to $BCDA$}
\FigureABCDequalsBCDA
\label{figure:ABCDequalsBCDA}
\end{figure}

\noindent{\em Proof}.  Refer to Fig.~\ref{figure:ABCDequalsBCDA}.  
To test whether $ABCD$ is equal to $BCDA$,  we place $AB$ horizontal
and make $BC^\prime = BC$ and $BA^\prime = AB$.  
Then all four of the smaller rectangles in the figure are congruent
(in the orientations shown).  Hence their diagonals are congruent;
hence $B$ is the common midpoint of the diagonals of the large
rectangle $DKD^\prime H$.  Hence $B$ is between $H$ and $K$.
Hence $ABCD$ is equal to $BCDA$.  That completes the proof.

\begin{lemma} \label{lemma:ER2} Two rectangles with the same height
are equal if and only if they have the same width.  Two rectangles 
with the same width are equal if and only if they have the same 
height.
\end{lemma}

\noindent{\em Proof}.   First observe that $a:b = a:c$ if and only 
if $b=c$,  as follows from the definition of proportionality.   Then 
the theorem follows from Lemma~\ref{lemma:ERproportion}.
That completes the proof.

\begin{lemma} \label{lemma:ER3} If one rectangle has both width 
and height less than a second rectangle, the two rectangles are 
not equal.
\end{lemma}

\noindent{\em Proof}. By Lemma~\ref{lemma:ERproportion},
it suffices to show that if $b<c$ and $a < d$,
then we do not have $b:c = d:a$.  To prove this, 
let $AOB$ be a right angle with $AO = c$ and $BO = d$.
By the definition of $<$ for line segments, there are points
$P$ and $Q$ with $\B(0,P,A)$ and $\B(0,Q,B)$ and 
$OP = b$ and $OQ = a$.  By inner Pasch, the lines
$PB$ and $AQ$ meet, and hence they are not parallel.
Then by definition, we do not have $b:c = d:a$.
That completes the proof.

\begin{lemma}\label{lemma:ER4} If equal rectangles are cut off
from equal rectangles (each cut off by one dividing line) then 
the remaining rectangles are equal.
\end{lemma}

\noindent
{\em Remark}. Euclid would have justified this by his common notion,
``if equals be subtracted from equals, the remainders are equal.''
\medskip

\begin{figure}[ht]
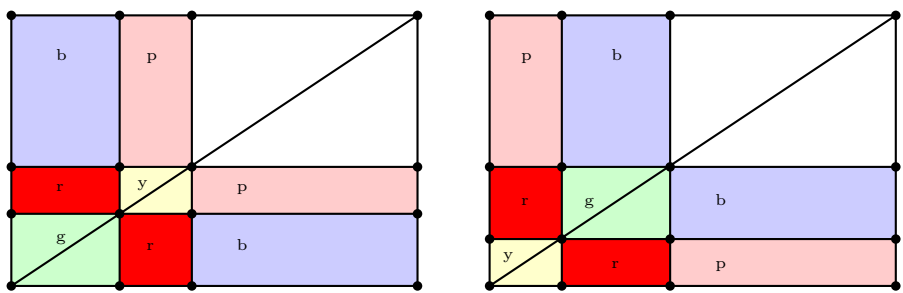

\caption{Cutting off equal rectangles (pink)
 from equal rectangles (pink and blue) leaves equal rectangles (blue).
 }
\label{figure:ER4}
\FigureERFourA
\FigureERFourB
\end{figure}

\noindent{\em Proof}.  
Refer to Fig.~\ref{figure:ER4}.  Rectangles
in the figure are labeled {\tiny p} for pink, {\tiny b} for blue,
{\tiny y} for yellow, and {\tiny r} for red, 
in case colors are not visible in the copy you are reading.
The hypothesis is illustrated
in the left part of the figure, 
 where the 
equal rectangles being cut off are pink, and the 
equal rectangles from which they are cut off are 
blue-and-pink.  Because these rectangles are equal,
the points shown on the diagonal are in fact on the diagonal.
Now we have to show that the blue rectangles are 
also equal.  By Lemma~\ref{lemma:ER1}, 
then pink-and-blue rectangles in the second figure 
are equal, so their common vertex lies on the diagonal.
The two green rectangles
are congruent, since their widths and heights are
equal to the width and height of the blue rectangles.
Similarly, the two yellow rectangles are congruent.
We have to show that the common point of the yellow and green
rectangles lies on the diagonal.  One easy proof of that is 
this: the left figure shows that the two red rectangles are equal.
But after permuting their vertices, those red rectangles are equal
to the the corresponding red rectangles in the right figure.
Therefore, the common point of the yellow and green rectangles 
lies on the diagonal.  That completes the proof.
\FloatBarrier

\begin{lemma}\label{lemma:ER5} If equal rectangles are pasted on to
equal rectangles, making larger rectangles, then 
the resulting larger rectangles are equal.
\end{lemma}

{\em Remark}. Euclid would have justified this by his common notion,
``if equals be added to equals, the wholes are equal.''  This can 
be taken literally if we define the sum of two rectangles having a common 
side to be their union. 
\medskip

\noindent{\em Proof.}  The theory of 
proportion will not be used.  We refer again to Fig.~\ref{figure:ER4}.
The rectangles assumed equal are the pink rectangles.  The 
rectangles pasted on are the blue rectangles. The lower left 
corner of the rectangle in the left figure is determined by the 
blue rectangles.  We must show that if the blue rectangles are equal,
that point lies on the diagonal of the white (and yellow) rectangles.
The equality of the blue rectangles is witnessed in the right figure
by the fact that the white and green rectangles have a common diagonal.
The green rectangles are congruent, and the yellow rectangles are congruent,
and the white rectangles are congruent. From the left figure, the 
yellow and white triangles have corresponding angles equal. 
From the right figure,
the green and white triangles have corresponding angles equal. 
Hence the yellow and green triangles have corresponding angles equal.  
Hence the yellow and green rectangles have 
a common diagonal (in both figures).  That completes the proof.

\subsection{Equal triangles}

\begin{definition}
The {\bf first circumscribed rectangle} $ABDK$ of triangle $ABC$
is the rectangle such that $C$ lies on the line through $D$ and $K$.
\end{definition}
\smallskip

(The word ``circumscribed'' in this context does not imply that the triangle
lies within the rectangle.)
To justify the definition, we have to show how to construct the 
circumscribed rectangle.  That is done as follows:  By I.29
there is a parallel line $L$ to $AB$ through $C$; then by I.12 there are  
perpendiculars $BD$ from $B$ to $L$ and $AK$ from $A$ to $L$.  Then
$AK$ is parallel to $BD$, since they are both perpendicular to $L$,
and $DK$ is parallel to $AB$ by construction, so $ABDK$ is a parallelogram.
Since $ABDK$ has a right angle, it is a rectangle.

\begin{figure}[ht]
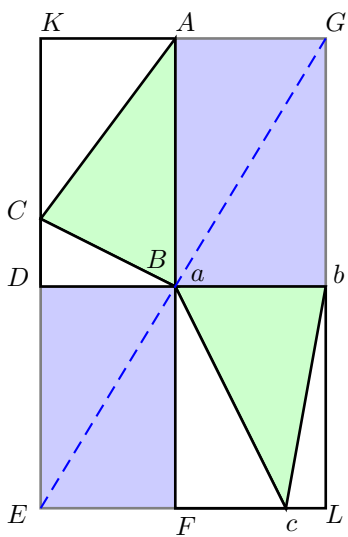

\caption{ $ABC$ is equal to $abc$ if $EBG$ is a straight line}
\label{figure:equaltriangles}
\begin{center}
\FigureEqualTriangles
\end{center}
\end{figure}

\begin{definition} Triangles $ABC$ and $abc$ are {\bf equal triangles} 
if their first circumscribed rectangles are equal.  See Fig.~\ref{figure:equaltriangles}.
\end{definition}

In this definition, the order of the vertices of $ABC$ matters.  
  To be precise (as we must when formalizing!) a triangle is 
an ordered triple of points.  Lemma {\tt ETpermutation} says that triangles
with the same vertices (in any order) are equal.  We will prove that below.

\begin{lemma} \label{lemma:ETequivalence} Equal triangles is an equivalence
relation.
\end{lemma}

\noindent{\em Proof}. The required properties of equal-triangles follow
from the definition and the fact that equal-rectangles is an equivalence
relation (Lemma~\ref{lemma:ERequivalence}).

\begin{lemma} [Euclid's I.37] \label{lemma:I.37}
Triangles which are on the same base and in the same parallels are 
equal to one another.  Specifically, they are equal triangles.
\end{lemma}

\noindent{\em Remarks}.  Remember that the base of $ABC$ is $AB$.
The base depends on the order in which the vertices are listed.
When applied to a segment,  the word ``equal'' means the same as 
``congruent'', and there is an $PQ = QP$, expressing the idea that
order is not important. Strict equality (identity) is never used,
so it is safe to use ``equal'' for ``congruent'', as Euclid does.
\medskip

\noindent{\em Proof}.  The base and ``parallels'' determine the circumscribed
rectangle up to congruence, so two triangles with equal (congruent) base and altitude have congruent
 circumscribed rectangles,  and congruent rectangles are equal. 
That completes the proof.

\section{Verification of the equal-triangles axioms}

\subsection{Remarks about the principle that $ABC$ is equal to $BCA$}
\label{section:VI.2}
This is one of the equal-figure axioms.
Euclid assumed this principle without proof (for the first time) 
at the penultimate line of the proof of I.35, where triangle $DGE$  has to be 
equal to triangle $GED$,  in order for the justification Euclid gives for 
that line (subtracting equals from equals) to apply. 

The principle is used again, and again without explicit mention,
in Euclid's Prop.~VI.2.  That proposition is equivalent to 
Corollary~\ref{lemma:proportion1}, so it shows that the principle 
that $ABC$ is equal to $BCA$ not only is implied by, but implies,
the basic theory of similar triangles.  Specifically, for those 
with their copy of Euclid at hand, in the proof of VI.2, triangle $BDE$
is first shown equal to $CDE$ since both have base $DE$ and the same
height,  and then is considered as having base $BD$ in order to 
apply V.1.  In the proof of V.1, applied here, $BDE$ is proved
equal to other triangles with base $BD$ and the same height.%

The referee offered to prove the equality of 
$ABC$ and $BAC$ as follows:  $ABC$ and $BAC$
are congruent by SAS (that is, I.4); and by one of our equal-figures 
axioms, congruent triangles are equal.  But we have 
already discussed in \S \ref{section:orientation} that ``order matters'',
and non-isosceles $ABC$ is not congruent to $BAC$.  
For this proof to be valid,  triangle congruence would have to be 
defined more generally, allowing the vertices to be mentioned in any 
order.  But that is neither the usage in Euclid, nor the tradition
in the intervening millenia, nor the modern formal usage.
So if we want $EBG$ and $BEG$ to be equal triangles in Euclid's proof
of I.35, we will have to justify that step somehow.  In
\cite{beeson2019} we introduced an axiom to do that.%
\footnote{There is a related precise question:  Is the 
equal-figure axiom that $ABC$ and $BAC$ are equal figures redundant?
That is, can it be derived from the rest of the axioms in the system 
of \cite{beeson2019}?  We do not know the answer.
} 

In this paper, we take a different approach: we 
supply a definition of ``equal triangles''.  Then we must
prove that $ABC$ is equal to $BCA$.  We
evidently cannot use Euclid's own theory of similar triangles
to justify the proof that $ABC$ is equal to $BCA$, since that would 
make the proof of VI.2 circular.  But we never intended to use 
Euclid's theory of similarity anyway,  as we wish to avoid the reliance
on Archimedes's axiom that is introduced in Book~V. 
We will instead use the nineteenth-century theory of 
proportionality.

Turning to modern times, the principle that $ABC$ is equal to $BCA$ is
also implicitly assumed and used in textbooks, as I will now explain. 
Consider the usual formula for calculating the 
area of a triangle:  base times height divided by 2.  Here we arbitrarily 
select one side as the base.  But if we select another side, will we necessarily
get the same answer?  ``Of course we will'', because the area depends only 
on the {\em unordered} vertices.  But in stating that, we have imported
some knowledge that is not in Euclid's axioms or the school curriculum.
If we were to try to {\em define} the 
area by the base times height over 2 formula,  
we would then have to {\em prove} that we get the same
answer no matter which side is taken as the base.  That is closely related
to the problem at hand, of verifying that $ABC$ is equal to $BCA$.  It is 
not exactly the same problem, as area is not involved in our definition 
of equal figures, but it addresses the same underlying issue.%

\subsection{Verification that $ABC$ and $BCA$ are equal}

\begin{theorem} \label{theorem:ETforward}
The  triangles $ABC$ and $BCA$ are equal.
\end{theorem}

\begin{figure}[ht]
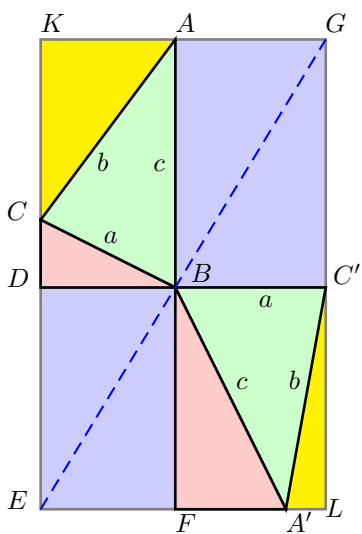

\caption{ $ABC$ is equal to $BCA$ if $EBG$ is a straight line}
\label{figure:ETpermutation}
\begin{center}
\FigureETpermutation
\end{center}
\end{figure}

\noindent{\em Proof}.  Applying the definition of equal 
triangles, we  get the situation
shown in Fig.~\ref{figure:ETpermutation}.

In the figure, the upper left and lower right rectangles are the 
circumscribed rectangles of the triangles, and the light blue rectangles
are constructed from them.  Their diagonals are $BG$ and $EB$.  It 
has to be proved that these diagonals lie on one straight line $EBG$.

We have
\begin{eqnarray*}
BA^\prime &=& c \ = \ AB \\
BC &=& a \ = \ BC^\prime \\
AB &=& b \ = \ AC^\prime \\
ABC &=& A^\prime B C^\prime \mbox{\qquad definition of triangle congruence} \\
\angle ABC &=& \angle A^\prime B C^\prime \ = \ \beta  \mbox{\ (the angle opposite $b$)}\\
\angle CBD &=& \angle A^\prime B F \ = \ \pi/2 - \beta 
\end{eqnarray*}

 Triangles $CDB$ and $A^\prime FB$ 
(shown pink in the figure, if the copy you are reading is in color) 
are similar, since they are both right 
triangles and have equal angles at $B$, as just shown. Therefore
$BF : BA^\prime = BD: BC$, or $BF:c = BD:a$.  
Since $BD = EF$, also $BF: c = EF :a$.
Then by Theorem~\ref{theorem:proportion-interchange},
$$BF:EF = c:a.$$
But $GC^\prime = AB = c $, and $BC^\prime = a$, 
so  $GC^\prime : BC^\prime = c: a$.    
Hence $$FB:EF = GC^\prime : BC^\prime,$$ 
by Lemma~\ref{lemma:proportion-transitive}.
Then by Corollary~\ref{lemma:proportional-implies-similar}, triangles $EBF$ and $BGC^\prime$ 
are similar.   Note that these triangles are right triangles, so that Corollary
 is applicable.  Hence 
$ \angle GBC^\prime =  \angle BEF$.  Since $EF$ is parallel to $DC^\prime$,
$EB$ extended past $B$ makes the same angle with $DC^\prime$ as with $EF$,
and hence coincides with $BG$.  Hence $EBG$ is a straight line. 
That completes the proof.

\subsection{Verification of the other equal-triangles axioms}

\begin{lemma}\label{lemma:congruentequal} Congruent
triangles are equal.
\end{lemma}

\noindent{\em Proof}.  Immediate from the definition of equal 
triangles, which is ``up to congruence.''
\medskip

\begin{lemma} \label{lemma:ETpermutation} Triangle $ABC$ is 
equal to triangles $BCA, CAB, ACB, BAC$, and $CBA$.
\end{lemma} 

\noindent{\em Proof}. We first prove that
the triangles $ABC$ and $BAC$ are equal.
The first circumscribed rectangles of 
$ABC$ and $BAC$ are congruent, since the bases $AB$ and $BA$ 
are equal, and the two altitudes are equal to the altitudes of 
triangles $ABC$ and $BAC$, which are equal.
Therefore $ABC$ and $BAC$ are equal, as claimed.
 
 Since all the permutations on three letters
are generated by $(213)$ and $(231)$, the conclusion follows from 
that case (which corresponds to permutation $(213)$) and 
Theorem~\ref{theorem:ETforward} (which corresponds to $(231)$
together with the transitivity of equal-triangles (Lemma~\ref{lemma:ETequivalence}).  That completes the proof.

The only other equal-figure axioms that involve only triangles 
are the two de Zolt axioms.  (These are the only axioms that assert
non-equality.  Without them, we could interpret all triangles as equal!)
The next two lemmas prove the de Zolt axioms as theorems, using the defined notion of equal triangles.

\begin{lemma} \label{lemma:deZolt1} 
Suppose $\B(b,e,d)$.  Then $dbc$ and $ebc$ 
are not equal triangles.  
\end{lemma}

\noindent{\em Proof}.  Suppose, for proof by contradiction,
that they are equal.  Then by the definition of equal triangles,
the first circumscribed rectangles of $dbc$ and $ebc$ are equal.
These rectangles have the same height (the perpendicular from 
$c$ to the line containing $d,b$, and $c$), but their widths
are respectively $db$ and $eb$.  By definition of less than for lines,
$ed < bd$.  But $bd = db$.  Since ``the part is not equal to the whole'',
this is impossible.  In our formal development, Euclid's principle 
that the part is not equal to the whole is a theorem about congruence.
By Lemma~\ref{lemma:ER2},
two rectangles of the same height are equal if and only if they 
have the same width; hence we have reached a contradiction.
That completes the proof.

\begin{lemma} \label{lemma:deZolt2} Let $abc$ be a triangle,
and suppose $\B(b,e,a)$ and $\B(b,f,c)$.  Then $abc$ is not 
equal to $ebf$.
\end{lemma}

\noindent{\em Proof}.  It suffices to show that triangle $bfe$
is not equal to $bca$.  The first circumscribed rectangles 
of $bfe$ and $bca$ have bases respectively $bf$ and $bc$, 
and $bf < bc$.  The altitude of $ebf$ is also less than the 
altitude of $abc$.  Hence, by Lemma~\ref{lemma:ER3}, they two 
rectangles are not equal.  Hence, by definition of equal triangles,
the triangle $bfe$ is not equal to $bca$. That completes the proof.
 
\section{Equal quadrilaterals}
Equal quadrilaterals will be defined in this section.
  The notion will be 
defined only for convex quadrilaterals and quadrilaterals
that are really triangles.    We will define the notion of 
convex so that it implies the quadrilateral lies in a plane; that is 
convenient, since we wish to do plane geometry, but without a 
dimension axiom,  following Euclid.%

The restriction to convex quadrilaterals bears some discussion.
There are two ways a quadrilateral might fail to be convex: either
it is {\em strictly} non-convex, or possibly it is ``really a triangle'',
i.e., one of its vertices is between the two adjacent vertices.
It seems that we must consider quadrilaterals that are really triangles,
as otherwise there is no connection between ``equal triangles'' and 
``equal quadrilaterals.''  That connection occurs in our former axiom 
(now to be a theorem) 
{\tt paste4} (which we will explain in due course, but not here). 
That axiom in turn is used in just one place:  Prop.~I.45.

In Euclid Books~I-III, 
 one does not find any propositions that mention or require
strictly non-convex quadrilaterals. Indeed if Euclid had wanted 
to include them, he would have had to be much more forthcoming
and explicit about what counts as a ``part'' of a figure, as a line
connecting two non-adjacent vertices can no longer be supposed to 
cut off a part.  In the rest of the paper, ``quadrilateral'' means
``convex quadrilateral'' or ``really a triangle.''

\begin{definition} Quadrilateral $JKLM$ is {\bf convex} if its 
diagonals meet.  That is, there is a point $E$ between $J$ and $L$
that is also between $M$ and $K$.
\end{definition}

\begin{definition} Let $JKLM$ be a convex quadrilateral.
 A {\bf circumscribed rectangle} of quadrilateral
$JKLM$ is a rectangle $ABCD$ with two sides parallel to 
diagonal $KM$, or two sides parallel to diagonal $JL$,
and each vertex of $JKLM$ lies on some side of $ABCD$.
\smallskip

A {\bf circumscribed rectangle} of a quadrilateral that is 
really a triangle is a circumscribed rectangle of that triangle.
\end{definition}

Each convex rectangle thus has two circumscribed rectangles, with 
sides parallel to one or the other diagonal of the quadrilateral.
A rectangle that is really a triangle has three circumscribed 
rectangles, all of which are equal rectangles by 
Theorem~\ref{theorem:ETforward}. 

\begin{definition}\label{definition:equalquadrilaterals}
 Two quadrilaterals are {\bf equal quadrilaterals}
or {\bf equal figures} if one of the circumscribed rectangles of one 
is equal to one of the circumscribed rectangles of the other.
\end{definition}

We turn now to the verification of the axioms for equal 
quadrilaterals.  See the Appendix for formal statements of 
these axioms. 

\subsection{Verification of the axioms {\tt EFPermutation} and 
        {\tt EFsymmetric}}

The axiom {\tt EFpermutation}  says that a quadrilateral is equal to any 
quadrilateral obtained by a permutation of the vertices, provided the 
two quadrilaterals have the same (set of) diagonals.
(That is, you can't just switch two adjacent vertices.)
  Since all such 
permuted quadrilaterals have the same two circumscribed rectangles, the 
verification of the axiom is immediate. 

The axiom {\tt EFsymmetric} follows immediately from the symmetry of the
relation of ``equal rectangles.''   The axiom {\tt EFtransitive} similarly
follows from the transitivity of ``equal rectangles''. 

\subsection{Verification of the axiom  {\tt halvesofequals}}
This axiom says 
that if equal quadrilaterals are each divided along a diagonal into 
equal triangles,  then all four triangles are equal.  Euclid used 
this principle without stating it%
\footnote{In Propositions~I.37, I.38, and I.48}
, and none of 
his common notions seem relevant.  

\begin{lemma} \label{lemma:halvesofrectangles}
The line connecting midpoints of opposite sides of a rectangle
divides it into two equal rectangles.
\end{lemma} 

\begin{figure}[ht]
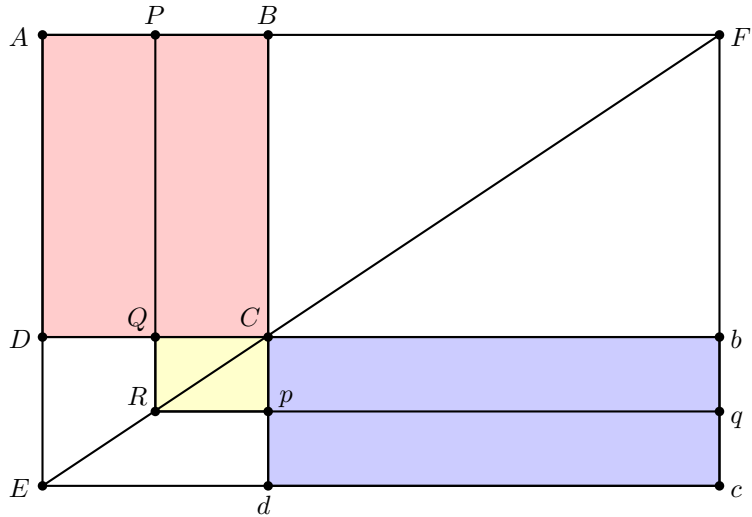

\caption{Equal rectangles $ABCD$ and $Cbcd$ are divided in half by $PQ$ and $pq$.  The resulting half-rectangles are equal since $R$ lies on the diagonal.}
\label{figure:halvesofrectangles}
\FigureHalvesOfRectangles
\end{figure}

\noindent{\em Proof}.  Let $ABCD$ and $abcd$ be the equal rectangles,
placed as shown in Fig.~\ref{figure:halvesofrectangles}, so that 
$C$ and $a$ coincide and the lines dividing the two rectangles in 
half are $PQ$ and $pq$ as shown.   Then since $ABCD$ and $abcd$
are equal rectangles, the sides can be extended to form a larger
rectangle, as shown, whose diagonal $EF$
passes through the point that is 
both $C$ and $a$.  Now let $R$ be the point of intersection of $PQ$
and $pq$.  Then $R$ is the intersection of the lines connecting
the midpoints of opposite sides of rectangle $DCdE$.  Then it is 
also the intersection of the diagonals of $DCdE$.  Then $R$
lies on the diagonal $EC$.  But that is collinear with $F$, since 
rectangles $ABCD$ and $CbcD$ are equal. Then by definition of 
equal rectangles, rectangles $PBCQ$ and $Cbqp$ are equal.
That completes the proof of the lemma.

\begin{lemma} \label{lemma:halveshelper}
 If a quadrilateral  $SPTQ$ is 
divided along its diagonal $PQ$ into two equal triangles, then 
the circumscribed rectangle $ABCD$ with $AD$ parallel to $PQ$
is divided into equal rectangles by $PQ$.  These rectangles
are also congruent.  
\end{lemma}

\noindent{\em Proof.}  By hypothesis, the triangles $SPQ$ and 
$TPQ$ are equal. (Since we have already verified axiom {\tt ETpermutation},
the order of the vertices does not matter.) Then by definition,
their circumscribed rectangles $APQB$ and $CPQD$ are equal.
Since they have the same height $PQ$, they are congruent by 
Lemma~\ref{lemma:ER2}. That completes the proof.

\begin{lemma}[Halves of equals are equal]\label{lemma:halvesofequals} If equal quadrilaterals
are each divided along a diagonal into 
equal triangles,  then all four triangles are equal.
\end{lemma}
\begin{figure}[ht]
\caption{The two quadrilaterals are equal, and 
each is divided in two equal triangles.  Then all the triangles are equal.}
\label{figure:halvesofequals}
\FigureHalvesOfEquals
\end{figure}

\noindent{\em Proof}.  See Fig.~\ref{figure:halvesofequals}.
 By Lemma~\ref{lemma:halveshelper},
the circumscribed rectangles of the two quadrilaterals are 
divided into two equal rectangles by the diagonals that 
divide the quadrilaterals into two equal triangles.  By 
Lemma~\ref{lemma:halvesofrectangles}, the half-rectangles are
all equal.  Then, by definition of equal triangles, the 
triangles whose circumscribed rectangles are those equal 
half-rectangles are equal. 
That completes the proof.
\FloatBarrier

\subsection{Verification of {\tt paste3} and {\tt paste4}}
\begin{figure} [ht]
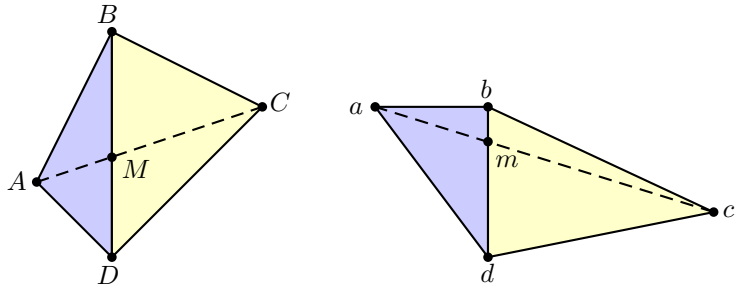

\caption{{\tt paste3}: If $ABD$ and $abd$   are equal, and $CBD$ and $cbd$   
are equal, then $ABCD$ and $abcd$ are equal.}
\label{figure:paste3}
\FigurePasteThree
\end{figure}

The axiom {\tt paste3} is the case of ``if equals be added to 
equals, the wholes are equals'',  when the equals being added
are triangles with a common edge, and the wholes are quadrilaterals.
See Fig.~\ref{figure:paste3}. 

\begin{figure} [ht]
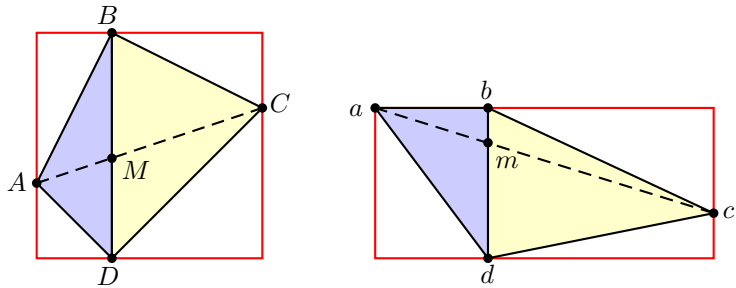

\caption{{\tt paste3}: showing the circumscribing rectangles}
\label{figure:paste3verification}
\FigurePasteThreeVerify
\end{figure}

\begin{lemma}[{\tt paste3}]\label{lemma:paste3} Suppose $ABC$ and $abc$ are
equal triangles, and $ABD$ and $abd$ are equal triangles, and 
$AC$ meets $BD$ at $M$ and $ac$ meets $bd$ at $m$.  The cases when 
$M=A$ or $M=C$ or $m=a$ or $m=c$ are also allowed. Then 
$ABCD$ and $abcd$ are equal quadrilaterals.
\end{lemma}

\noindent{\em Proof}.  Construct the circumscribed rectangles of 
$ABCD$ and $abcd$ (with sides parallel to $BD$ and $bd$, respectively).
 Then these rectangles are each divided into 
two rectangles, which are the circumscribed rectangles of triangles
$ABD$, $BCD$, $abd$, and $bcd$.  Since the circumscribed rectangles
of equal triangles are equal, the hypotheses of Lemma~\ref{lemma:ER5}
are satisfied, and the conclusion of that lemma is that
the  circumscribed rectangles of $ABCD$ and $abcd$ are equal rectangles.
Then by Definition~\ref{definition:equalquadrilaterals}, $ABCD$
is equal to $abcd$.  That completes the proof.
\smallskip

{\em Remark.} The axiom {\tt paste3} was used in \cite{beeson2019} in the proofs of 
Propositions~I.35a, I.42, I.47B, and in proving lemma {\tt paste3},
which was in turn used in I.48.  In the present approach, we prove
I.42 using {\tt doublesofequals} (which will be proved next)
and the rest directly.  However, now we use {\tt paste3} in 
proving I.45.  

\begin{lemma}[Doubles of equals are equal]\label{lemma:doublesofequals}
If two quadrilaterals $ABCD$ and $abcd$ are each divided along the diagonals
$AB$ and $ab$ into two
equal triangles, and  triangles $ABC$ and $abc$ are
equal, then the 
two quadrilaterals are equal.
\end{lemma}

\noindent{\em Remark}. It is allowed that one or both quadrilaterals
may be ``really a triangle'', with one end of the diagonal between the 
adjacent vertices.  This is needed in proving I.44.
\medskip

\noindent{\em Proof}.  This lemma is immediate consequence of
{\tt paste3}. However, we like the following proof directly from the definition,
that does not appeal to Lemma~\ref{lemma:ER5}, so we present it as well.
 Refer to Fig.~\ref{figure:halvesofequals}. 
The hypothesis is that all four triangles shown are equal. 
We have to show that rectangles $ABC$ and $Cbcd$ are equal, which 
implies that the red quadrilateral and the green quadrilateral are equal.
Because the red and green triangles are equal, point $R$ lies
on the (extended) diagonal $FC$. Because the two red triangles are equal,
$Q$ is the midpoint of $DC$.  Because the green triangles are equal,
$p$ is the midpoint of $Cd$.  The diagonals of a rectangle bisect each 
other (as Euclid could have easily proved after Prop.~I.34 without using any
equal figures axioms; we did so in \cite{beeson2019}); the diagonals of 
$DCdE$ in particular meet at $R$.  Hence $C$ lies on $EF$.  Hence 
rectangle $ABCD$ is equal to rectangle $Cbcd$, as claimed.  That 
completes the proof.
\smallskip

The axiom {\tt paste4} is very similar to {\tt paste3}:  like {\tt paste3}, it is about pasting
two triangles together to get a quadrilateral.  But in {\tt paste4}, the ``triangles''
are quadrilaterals that are ``really triangles'', in the sense that one vertex is 
between the two adjacent vertices.   The picture is the same as Fig.~\ref{figure:paste3},
except that two more points are added along two sides of the triangles.  The verification
that {\tt paste4} holds with the defined notion of equal quadrilaterals is the same 
as the verification of {\tt paste3},
since the circumscribed rectangle of such a quadrilateral is 
the same as the circumscribed rectangle of the triangle.
\smallskip

{\em Remarks}.
In our formalization \cite{beeson2019},  Axiom 
{\tt paste3} was used to prove Propositions I.35 and Lemma
{\tt EFreflexive}, both of which 
we can now prove directly, and to prove Propositions~I.42 and I.42B
and I.47B
where {\tt doublesofequals} suffices.  In fact in the proof of I.47,
Euclid specifically states, in parentheses, ``But the doubles of equals
are equal to one another.''  This is clearly intended as a justification
for the statement that follows that quotation, which however has no official justification,
as Euclid has not proved that doubles of equals are equals.

\subsection{Proposition I.35}
\begin{lemma}[Proposition I.35] \label{lemma:I.35}
Parallelgrams which are on the same base and in the same parallels
are equal to one another
\end{lemma}

\begin{figure}[ht]
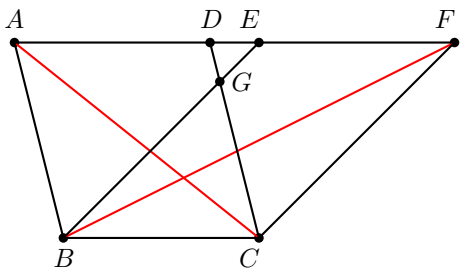

\caption{Euclid's figure for I.35}
\label{figure:I.35}
\FigureOneThirtyFiveWithDiagonals
\end{figure}
\medskip

\noindent{\em Remark}. This is the first proposition in which 
Euclid needs (something like) the equal-figures axioms; in particular
triangle $DEG$  (see Fig.~\ref{figure:I.35}) 
needs to be equal to triangle $EDG$ for Euclid's proof
to work.  But here we derive Prop.~I.35 from I.37, which we have derived
already, directly from the definition of ``equal triangles.''
\medskip

\noindent{\em Proof}.  Fig.~\ref{figure:I.35} is Euclid's figure, except that
we have added the (red) diagonals $AC$ and $BF$.
By I.37, triangles $ABC$ and $BCF$ are
equal, since they are triangles with the same base and in the 
same parallels. By I.34, the diagonals divide the parallelograms
into equal triangles.  By Lemma~\ref{lemma:doublesofequals}, the
parallelograms are equal.  That completes the proof.
\FloatBarrier

\subsection{Prop.~I.42, I.43, and I.44}
We have already remarked that I.42 follows from Lemma~\ref{lemma:doublesofequals}.  I.44 does not use any 
equal-figures axioms, other than transitivity, which we have proved.
Therefore I.44 is provable. 
  We will turn to I.43, but we 
need a lemma first.

\begin{lemma} \label{lemma:equalrectanglesequalfigures}
Equal rectangles (in the sense of Definition~\ref{definition:equalrectangles})
are equal quadrilaterals (in the sense of Definition~\ref{definition:equalquadrilaterals}).
In other words, circumscribed rectangles of equal rectangles are equal. 
\end{lemma}

\begin{figure}[ht]
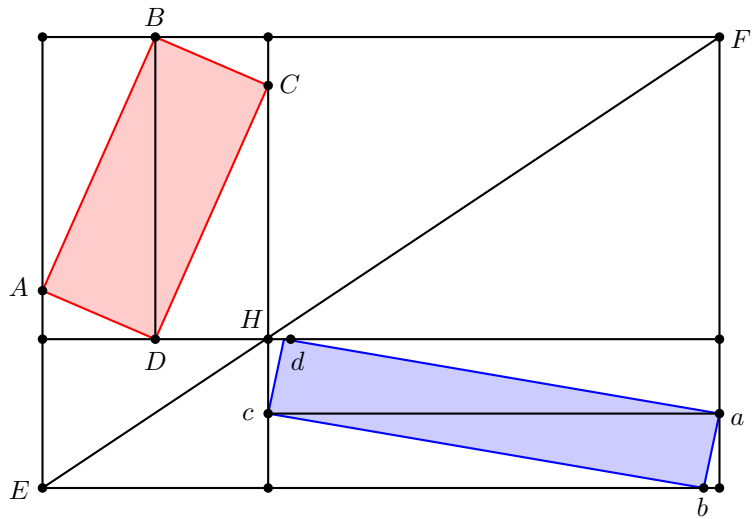

\caption{Equal rectangles are also equal quadrilaterals}
\FigureEqualRQ
\end{figure}
\medskip

\noindent{\em Proof}.  Let $ABCD$ and $abcd$ be equal rectangles.
Their diagonals divide each of them into two congruent triangles.
By the definition of equal triangles,  triangle $ABD$ is equal to 
triangle $abc$, since they have equal circumscribing rectangles 
(and we have proved the order of vertices does not matter). 

 One circumscribed rectangle of triangle $ABD$ is half of the circumscribed rectangle 
of $ABCD$, and one circumscribed rectangle of triangle $abc$ is half
of the circumscribed rectangle of $abcd$.  
Then by Lemma~\ref{lemma:doublesofequals}, quadrilaterals $ABCD$ and $abcd$ 
are equal.  That completes the proof. 

\FloatBarrier
\begin{lemma}[Prop.~I.43] \label{lemma:I.43}
In any parallelogram the complements of the parallelograms about the diameter
are equal to one another.
\end{lemma}

\noindent{\em Remark}.  In Euclid's proof, he needs to apply the common 
notion that equals subtracted from equals leave equal remainders.
In formalizing that argument in \cite{beeson2019}, we used the 
axioms {\tt cutoff1} and {\tt cutoff2}.  But here, we prove the theorem
without common notions or cutoff lemmas, by reducing it to the case of 
a rectangle instead of a parallelogram, and then using the definition
of equal rectangles.
\medskip

\begin{figure}[ht]
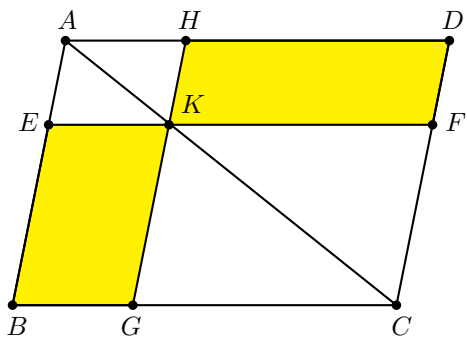

\caption{Euclid's figure for Prop.~I.43. To prove: the yellow 
parallelograms $BGKE$ and $KFDH$ are
equal.}
\label{figure:I.43}
\FigureOneFortyThree
\end{figure}

\noindent{\em Proof}. Let $A^\prime$ be a point collinear with $AD$ 
such that $A^\prime B$ is perpendicular to $BC$.  Let $E^\prime$ 
be the intersection of $A^\prime D$ and the line containing $EF$. 
Similarly let $H^\prime K^\prime$ be perpendicular to $BC$ at $G$ 
and $D^\prime F^\prime$ perpendicular to $BC$ at $C$, with 
$K^\prime F^\prime$ collinear with $KF$ and $H^\prime D^\prime$
collinear with $HD$.  Since opposite sides of a parallelogram are 
equal, and parallelograms with equal bases and in the same parallels
are equal, the parallelograms with primes are equal to the parallelograms
without primes.  Therefore, it suffices to prove the lemma in the 
case when all the parallelograms are rectangles. But in that case,
by the definition of equal rectangles, $H^\prime D^\prime K^\prime F^\prime$
is equal to $E^\prime K^\prime G^\prime B^\prime$, as rectangles,
if and only if $K^\prime$ lies on $A^\prime C$.  Now by Lemma~\ref{lemma:equalrectanglesequalfigures}, they are also equal in the 
sense of ``equal quadrilaterals.'' That completes the proof.

\subsection{Prop.~I.45}

In the formalization of \cite{beeson2019},  Prop.~I.45 is proved with the aid of the axiom {\tt paste2} and
the theorem {\tt paste5},  which was originally an axiom, but can be proved from the other paste axioms. 
The paste axioms are used to 
 formalize an inference that Euclid justifies 
by the common notion ``if equals be added to equals, the wholes are equal.''
Since we have not verified {\tt paste2},  we either have to do so, or prove I.45 directly.
In this section we give a direct proof of  the step of I.45 in question.

\begin{figure}[ht]
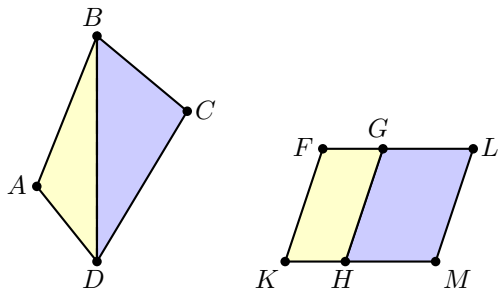

\caption{Euclid's figure for Prop. I.45 (color added)}
\label{figure:I.45}
\FigurePasteFiveHelperG
\FigureFLMKParallelogram
\end{figure}

In Euclid's proof of I.45, we have a quadrilateral $ABCD$,  whose two parts $ABC$ and $ABD$
are respectively equal to parallelograms $FGHK$ and $GLMH$, which share the 
common side $GH$ and together form the ``whole'', that is, the larger parallelogram 
$FLMH$, as shown in Fig.~\ref{figure:I.45}.

By Euclid~I.36, we could reduce to the case when the parallelograms are rectangles,
but that does not make the proof immediate.  Axiom {\tt paste2} permitted,
in our old formalization,  
``adding'' triangles that share a common side,  but here we have to ``add'' 
equals,  where on one side we are adding parallelograms (or rectangles) with a common side.

That is the purpose of the theorem  {\tt paste5}, which we proved in \cite{beeson2019} from the 
axiom {\tt paste2}.     More generally,   {\tt paste5} permits
 $FLMH$ to be
any convex quadrilateral, not just a parallelogram, divided in two by a line $GH$ from 
side $FL$ to side $KM$. In this section, we only prove the parallelogram case, which 
suffices for I.45. 

\begin{lemma}\label{lemma:paste5helper}
Every  quadrilateral  or triangle is equal to half its circumscribed rectangle,  where ``half'' 
means that the rectangle is divided by connecting the midpoints of two opposite sides.
\end{lemma}

\begin{figure}[ht]
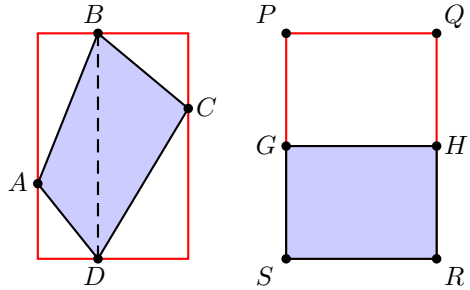

\caption{$ABCD$  is equal to half of $PQRS$, namely $GHRS$}
\label{figure:paste5-1}
\FigurePasteFiveHelperA
\FigurePasteFiveHelperB
\end{figure}

\begin{figure}[ht]
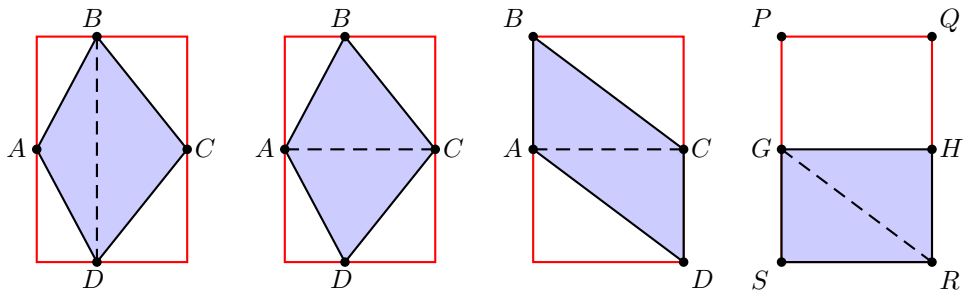

\caption{The blue (or shaded) quadrilaterals are all equal.}
\label{figure:paste5-2}
\FigurePasteFiveHelperC
\FigurePasteFiveHelperD
\FigurePasteFiveHelperE
\FigurePasteFiveHelperF
\end{figure}

\noindent{\em Proof}.   By Euclid~I.40 (triangles with the same base and in the same parallels are equal),
we may without 
loss of generality assume that triangles $BDA$ and $BDC$ are isosceles, i.e., 
$AB=AD$ and $BC = CD$.   Applying Euclid~I.40 again, we may assume without loss of 
generality that $B$ is in the upper left corner and $D$ in the lower right corner.  See Fig.~\ref{figure:paste5-2}.

Now we are in a position to apply {\tt paste3} (Lemma~\ref{lemma:paste3}).   Triangle $ACD$ 
is equal to triangle $GHS$ (since they are congruent).   Triangle $ACB$ is equal to triangle $GSR$
(since they too are congruent).   By Lemma~\ref{lemma:ETpermutation}, the order of the vertices does 
not matter.  By {\tt paste3}, quadrilaterals $ABCD$ and $GHSR$ are equal; again the order of the 
vertices does not matter, by the definition of equal quadrilaterals.   That completes the proof.

\begin{lemma}\label{lemma:addequals}
Let $ABCD$ be a convex quadrilateral and $FLMK$ a parallelogram,
composed of two parallelograms $FGHK$ and $GLMH$.  Suppose triangle
$ADB$ is equal to $FGHK$ and triangle $CBD$ is equal to $GLMH$.
Then $ABCD$ is equal to $FLMK$.  See Fig.~\ref{figure:I.45}.
\end{lemma}

\noindent{\em Remark}.  This is the part of Euclid~I.45 that 
is justified using the common notion ``if equals be added to equals,
the wholes are equal.''  In our formalization \cite{beeson2019},
we used an equal-figures axiom {\tt paste4}. 
\smallskip

\noindent{\em Technical remark}. We did not formally define the notion of 
a triangle being equal to a quadrilateral.  To state the lemma 
formally, we need to
consider another point $E$ between $B$ and $D$
somewhere on the boundary of each triangle mentioned, so triangles
$ABD$ and $BCD$ are formalized as quadrilaterals $ABED$ and $CBED$. 
\medskip

\begin{figure}[ht]
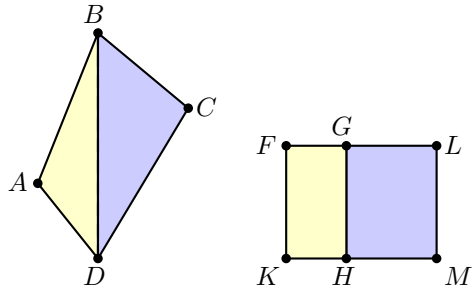

\caption{To prove: $ABCD$ is equal to $FLMK$.}
\label{figure:paste4}
\FigurePasteFiveHelperG
\FigureFLMKRectangle
\end{figure}

 \noindent{\em Proof}. By Euclid~I.36 (which we have already proved),
it suffices to consider the case in which the parallelograms are rectangles.  The 
resulting situation is shown in Fig.~\ref{figure:paste4}.

\begin{figure}[ht]
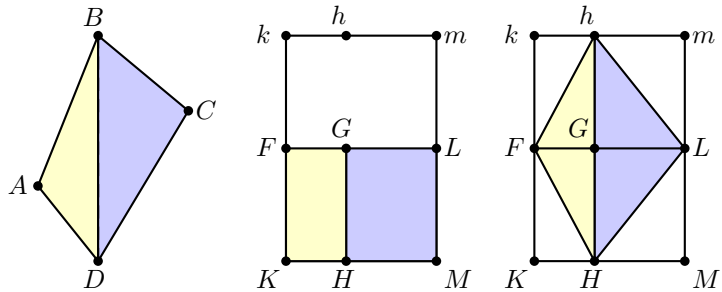

\caption{Some steps in the proof of Lemma~\ref{lemma:addequals}}
\label{figure:paste4b}
\FigurePasteFiveHelperG
\FigureFLMKBigRectangle
\FigureFLMKBigRectangleB
\end{figure}

We construct a larger rectangle, double the height of $FLMK$, by 
adding a new rectangle $kmLF$ above $FL$, the same size and shape as 
$FLMK$. Also define $h$ on line $GH$ and line $km$.   Then $G$ is the midpoint of $hH$
and $L$ is the midpoint of $mM$. 
See Fig.~\ref{figure:paste4b}.

Triangle $FGh$ is congruent to $HKF$ (since the 
diagonals of a parallelogram divide into two congruent triangles), and 
also to triangle $FGH$, since right triangles with congruent legs are
congruent.  Then by Lemma~\ref{lemma:paste3} ({\tt paste3}), 
$FGHK$ is equal to $hGHF$.  Since $hGHF$ is really a triangle, 
$FGHK$ is equal to $hHF$.   Therefore $ABD$ is equal to $hHF$.
Similarly, $ABD$ is equal to $hHL$.   Then by  Lemma~\ref{lemma:paste3} ({\tt paste3}), $ABCD$ is equal to 
$FhLH$.   By Lemma~\ref{lemma:paste5helper},  $FhLH$ is equal to half its circumscribed
rectangle, which is $FLMK$.   By the transitivity of ``equal quadrilaterals'', $ABCD$ is equal to $FLMH$.
That completes the proof.
\medskip

\begin{corollary} \label{lemma:paste4}  Euclid I.45 is  provable,
either abstractly from the other equal-figure axioms, but wihout {\tt paste2} or any of the 
cutoff axioms,   or using the definitions of 
``equal triangles'' and ``equal quadrilaterals'' and not the equal-figure axioms.
\end{corollary}

\noindent{\em Remark}.  This corollary eliminates the need to prove {\tt paste5} as part of the 
formal development of Euclid.   It does not quite eliminate {\tt paste2}, as that is used in 
the formal development to prove I.35  as well as {\tt paste5}.   Here we proved I.35 
directly,  using the defined notion of ``equal figures'',  so {\tt paste2} was not needed anywhere
in Euclid Book~I.
\medskip

\noindent{\em Proof}.   Lemma~\ref{lemma:addequals} provides the step in the proof of Prop.~I.45
for which we needed {\tt paste4}.   That completes the proof.
\medskip

\section{A Euclidean theory of area}

Prop.~I.45  allows to construct a parallelogram with  one angle specified that is equal to 
a given quadrilateral.   The same method of proof, taken one step further, could have allowed
Euclid to construct a parallelogram with one angle {\em and one side} specified,  equal to a given
quadrilateral.  Perhaps Euclid thought that was so obvious that he did not need to spell it out.  
That is somewhat surprising,  since it could be used to define ``equal area.''  Namely, two figures
have equal area, if they are both equivalent to the same rectangle.  If one side of the rectangle is 
regarded as a linear measuring unit, the other side gives the area of the figure measured in square units. 
But Euclid did not take this step.  We take it in the following definition:

\begin{definition} Let $\Delta$ be a convex quadrilateral or triangle,
and let $UV$ be any given segment (thought of as the ``unit segment''). 
Then the {\bf area of $\Delta$  relative to unit $UV$} is 
any rectangle equal to $\Delta$ and having   one side equal to $UV$. 
\end{definition}

\begin{theorem}\label{theorem:area} Every convex quadrilateral or 
triangle has an area, and the area is unique (up to congruence). 
\end{theorem}

\noindent{\em Proof}.  As remarked above, this is a minor 
extension of Euclid~I.45. 
\smallskip

 By the ``sum'' of two
rectangles with side $UV$ we mean another such rectangle formed by placing the two end-to-end; of course that
notion can be formulated precisely without mentioning motion.   With $UV$ fixed,  the ``area'' of a triangle or quadrilateral $\Delta$ is 
by definition the rectangle 
with side $UV$ that is equal to the circumscribed rectangle of $\Delta$.  That rectangle is unique, 
by Lemma~\ref{lemma:ER2}.
 We can now prove that we can test two figures for equality by comparing their areas: 

\begin{theorem} \label{theorem:area2}  With area defined as above,  relative to a fixed ``unit segment'',
two  quadrilaterals or triangles are equal if and only if  they have the same area.
\end{theorem}

\noindent{\em Proof}. Let $\Gamma$ and $\Delta$ be equal figures. 
As remarked above, the method of
 Prop.~I.45  and I.~38 can be used to show that each figure is equal to a rectangle with one side equal to the
 unit segment $UV$.   By transitivity,  the figures are equal if they have equal area.   That is half the theorem.
 Now, suppose the figures $\Gamma$ and $\Delta$ have equal area; we must prove they are equal.  Since they have equal area, 
 by the definition of area,  there are rectangles $R_1$ and $R_2$, equal  respectively to the circumscribed rectangles of $\Gamma$ and $\Delta$,
such that $R_1$ and $R_2$ are equal rectangles.    By Lemma~\ref{lemma:equalrectanglesequalfigures}, $R_1$ and $R_2$ are 
equal quadrilaterals (as well as equal rectangles).  By the transitivity of ``equal rectangles'',  the circumscribed rectangles of $\Gamma$ and
$\Delta$ are equal.  Hence, by definition, $\Gamma$ and $\Delta$ are equal figures.  That completes the proof.
 \medskip
 
 \noindent{\em Remark}.  This is a kind of ``completeness theorem'' for the notion of equal figures, and also for the
 equal-figures axioms.   If ``figures''  means convex quadrilaterals and triangles, then we have not omitted any  axioms that
 would be needed to use ``equal area'' as an interpretation of ``equal figures''.   But that is only one property we would wish 
 ``area'' to have. 
 \medskip
 
 \begin{theorem}[Additivity of area] \label{theorem:additivity}
 Suppose figures $\Gamma$ and $\Delta$  have a common line $BD$ and together form another figure $\Xi$. 
 Then the  sum of the areas of $\Gamma$ and $\Delta$ is the area of $\Xi$.
 \end{theorem}
 
 \noindent{\em Remarks}.  Remember that the area is a rectangle (with unit height) and  ``sum'' refers to the ``rectangle sum'',
 resulting from placing two unit-sided rectangles end-to-end with unit sides together.  Also remember that ``figure'' means (for the proof below at least) 
 triangle or convex quadrilateral.   This theorem amounts to the ``mother of all paste axioms,''  in the sense that it 
 represents more closely the intuition behind the paste axioms.
 \medskip
 
 \noindent{\em Proof}.    Let $R_1$ be the area of $\Gamma$ .   Let $R_2$ be the area of $\Delta$.    (Then $R_1$ and $R_2$ 
 are rectangles.)  Let us first consider the case when $\Gamma$ and $\Delta$ are both triangles (or quadrilaterals that are 
 ``really triangles'')  sharing an edge. 
 Since $\Xi$ is a figure, it is either a convex quadrilateral or ``really a triangle''.   Since $R_1$ is equal to $\Gamma$ 
 (by definition of area) and $R_2$ is equal to $\Delta$,  we can conclude by  Lemma~\ref{lemma:addequals} 
 that $\Xi$ is equal 
 to the rectangle sum $R_1 + R_2$, i.e., the rectangle obtained by pasting together $R_1$ and $R_2$ along their unit edges.
 But by Lemma~\ref{lemma:ER3}, that rectangle must be the area of $\Xi$.  That completes the proof in case $\Gamma$
 and $\Delta$ are triangles.  
 
 Before proceeding,  we note that rectangle addition $R_1 + R_2$ is both commutative and associative.  These
 properties follow from Lemma~\ref{lemma:ER2} and elementary theorems about congruence of line segments.

\begin{figure}[ht]
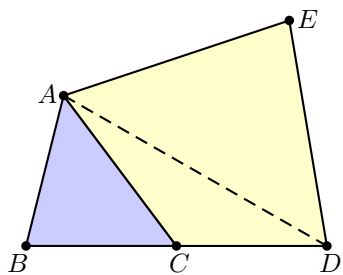

 \caption{$area(ABDE) = area(ABC) + area(ACDE)$}
 \label{figure:additive1}
 \FigureAdditiveOne
 \end{figure}

 Now consider the case when $\Gamma$ is a triangle and $\Delta$ is a convex quadrilateral.  Since $\Xi$ is a quadrilateral,
 $\Gamma$ and $\Delta$ must have a vertex in common.  The situation must be as shown in Fig.~\ref{figure:additive1},
  in which $\Xi$ is $ABDE$,  $\Gamma$ is $ABC$, and $\Delta$ is $ACDE$.     
 By hypothesis, $\Xi$ is a figure, so it is convex, i.e., $B$ meets $AD$.
Then line $AD$ divides $\Delta$ into 
 two triangles, and we have 
 \begin{eqnarray*}
 area(\Xi) &=& area(AED) + area(ADB) \\
 &=& area(AED) + (area(AEC)  + area(ACB)) \\
 &=& (area(AED) + area(AEC)) + area(ACB)  \mbox{\qquad see below} \\
 &=& area(AEDC) + area(ACB)  \mbox{\qquad by associative of rectangle addition} \\
 &=& area(\Delta) + area(\Gamma) \\
 &=& area(\Gamma) + area(\Delta)
 \end{eqnarray*}
 That completes the proof in case $\Gamma$ is a triangle and $\Delta$ is a convex quadrilateral.  Since rectangle
 addition is commutative, that also takes care of the case when $\Delta$ is a triangle.  
 
 The only remaining case is when both $\Gamma$ and $\Delta$ are quadrilaterals.   Then the situation is as 
 shown in Fig.~\ref{figure:additive2}. 
 
 \begin{figure}[ht]
 \caption{$area(ACDE) = area(ABEF) + area(BCDE)$}
 \label{figure:additive2}
 \FigureAdditiveTwo
 \end{figure}
We apply the same method:  adding a quadrilateral amounts to adding
 two triangles.   Then 
 \begin{eqnarray*}
 area(\Xi) &=& area(ACDF) \\
 &=& area(ACDE) + area(AFE)  \mbox{\qquad as proved above} \\
 &=& (area(ABE) + area(BCDE)) + area(AFE) \\
 &=& (area(ABE) + area(AFE)) + area(BCDE) \\
 &=& area(ABFE) + area(BCDE) \\
 &=& area(\Gamma) + area(\Delta)
 \end{eqnarray*}
That completes the proof.
\medskip

\noindent{\em Remark}.  It is worth remembering that our treatment has 
equal-figure axioms only for triangles and convex quadrilaterals. 
While that seems to cover Euclid Books I to II, it is still unsatisfactorily far from capturing the intuitive notion of ``figure'', which
should at least include all (convex or not) connected polygons, and perhaps disconnected ones,  and ones with curved boundaries.  In particular, 
Book~IV deals with polygons inscribed in circles, so it will be necessary
to extend the notion of ``figure'' before we can formalize Book~IV. 
These matters are discussed at length in \cite{hartshorne} and \cite{hilbert1899}, and are beyond the scope of this paper, whose
aim is simply to remove the necessity for an undefined notion of ``equal figures'' in Euclid Books I  and II (and III comes for free, since it 
does not mention figures at all).  
    
\subsection{A conservative extension theorem}
Let us review.  We have defined ``equal triangles'' and ``equal 
quadrilaterals'', and proved many of the equal-figure axioms, using
those defined notions instead of a primitive notion.   In particular we 
proved that those defined notions are equivalence relations, and that
the order of listing the vertices does not matter,  and we proved
{\tt paste3} and {\tt paste4},  corresponding to pasting triangles along a 
common side.   

We so far never needed the following axioms: {\tt paste1}, {\tt paste2}, {\tt cutoff1},
and {\tt cutoff2}.  In fact, as you can see in the appendix:
\begin{itemize}
\item  {\tt paste1} was never
used, even in the formalization in \cite{beeson2019}.
\item [\tt paste2]  was 
used (only) to prove Prop.~I.35,  which we proved directly in Lemma~\ref{lemma:I.35},
and to prove {\tt paste5}, which in turn was used (only) in Prop.~I.48.
\item{\tt cutoff1} and {\tt cutoff2} were used (only) in Prop.~I.35 and Prop.~I.43, 
both of which we have proved directly in this paper.
\end{itemize}

 We also were able
to prove directly those propositions of Euclid in which he made
use of ``equal figures.''  Therefore, we have now achieved the aim 
of showing that Euclid did not actually need to take the notion of 
``equal figures'' as primitive.%
\footnote{Except of course, that we still need to show how to develop
the theory of proportion with techniques from Book I.}

But there is still a bit more to prove concerning the relation between the defined 
notion of ``equal figures'' and the axiomatic treatment in \cite{beeson2019}.

\begin{theorem} [Conservative extension theorem] 
\label{theorem:conservative} 
Any theorem not mentioning ``equal figures'' that can be proved
with the aid of the equal-figures axioms,  can also be proved without 
those axioms.
\end{theorem}

\noindent{\em Proof}.   More generally,  in any formal proof we can replace
the symbols {\tt ET} and {\tt EF}  by the defined notions of ``equal triangles''
and ``equal quadrilaterals'',  to get the ``interpretation'' of the proof.  We 
claim that every step in the interpretation of a proof is a theorem provable 
without equal-figure axioms.   This we prove by induction on the length 
of proofs.   The only difficult part is the base case, in which $\phi$ is an 
equal-figures axiom.   In this paper we have proved the interpretations of 
all the equal-figures axioms except {\tt paste1} and {\tt paste2}.  We must
now take care of this loose end. 

Fortunately, we are spared from the need to prove {\tt paste1} and {\tt paste2}
directly from the definitions of ``equal triangle'' and ``equal quadrilateral''.
Instead, they can be derived from the additivity of area.   Here is how.
Both those axioms have the form that some figure $\Xi$ is 
composed of two figures $\Gamma$ and $\Delta$ where $\Gamma$ 
is a triangle and $\Delta$ is a quadrilateral.  In {\tt paste1}, the two 
are combined as in Fig.~\ref{figure:additive1}, and in {\tt paste2},
as in Fig.~\ref{figure:additive2}.   The hypothesis is that there
are figures $\Gamma^\prime$ and $\Delta^\prime$, equal respectively 
to $\Gamma$ and $\Delta$,  and combining to make a figure $\Xi^\prime$.
The desired conclusion is that $\Gamma^\prime$ and $\Delta^\prime$ 
are equal figures.  According to Theorem~\ref{theorem:area2},
it suffices to prove that $area(\Xi ) area(\Xi^\prime)$.
That follows from the additivity of area.  Explicitly, we have
\begin{eqnarray*}
area(\Xi) &=& area(\Gamma) + area(\Delta) \\
&=& area(\Gamma^\prime) + area(\Delta^\prime) \mbox{\qquad by Theorem~\ref{theorem:area}} \\
&=& area(\Xi^\prime) \mbox{\qquad by the additivity of area}
\end{eqnarray*}
That completes the proof of the theorem.

\section{The early theory of proportionality} \label{section:kupffer}

 The question of 
``early development of the theory of proportions'' (that is, 
 its development 
without the Axiom of Archimedes, from principles in the first part of 
Euclid Book~I) 
was already
considered by German geometers in the 1890s.%
\footnote{Thus ``early'' refers to the logical status of the work, i.e.,
not relying on the later parts of Euclid that need Archimedes, and not to 
the chronology, since this work was done in the nineteenth century.}  
 This subject has been
worked on by Kupffer, Schur, Dehn, Hessenberg, Hilbert, and Bernays. 
Above we showed that the heart of the matter is just two theorems, 
which we call (following Bernays) the interchange theorem and the 
fundamental theorem of proportion.  
In this section,  we discuss the beautiful proofs of Karl Kupffer, who found them in 
1893, as described in \cite{schur1902}, and those of Bernays, both
of which answer to our need for a development by the means available to Euclid.%
\footnote{Already in 1810, Bolzano \cite{bolzano1810} called 
Euclid's use of Book~V to reach the theory of similar triangles an 
``atrocious detour''.  Baldwin \cite{baldwin2018} calls this 
``Bolzano's challenge'', and compares the treatments of similarity in 
Euclid, Descartes, and Hilbert.} 

\subsection{Proportionality in Euclid~Book VI and in Hilbert}
Euclid Book~VI presents a theory of proportionality 
based on the theory of 
``magnitudes'',  due to Eudoxes and using the axiom of Archimedes, which
occupies Book~V.   
The interchange theorem does occur in Book~VI,  but  Book~VI makes use of Prop.~I.38, 
which makes use of 
I.~35, which requires the use of the equal-figures axiom that $ABC$
is equal to $BCA$.  It would therefore be circular to make I.35 rely on 
the Eudoxian theory of proportions.   Perhaps this is the reason that Euclid turned to 
applying the common notions to ``equal figures.''

The theory of proportion
can be derived from two fundamental theorems of projective geometry, 
Desargues's theorem and Pascal's theorem. 
Hilbert \cite{hilbert1899}, \S14,  sketches a proof of Pascal's theorem based on 
theorems about circles occurring in Euclid Book~III.  To use Hilbert's
approach, we would have to prove those theorems about circles and 
cyclic quadrilaterals, without using the equal-figure axioms of Euclid.
That is possible, but it would require inserting a few theorems 
from Book~III before Prop.~I.35.   However, it is not at all certain
that Desargues's theorem can be proved with the tools at hand 
before I.35.  Hilbert only proves it in \S24, after defining segment
arithmetic.  Thus, we do not find direct support in Hilbert for the claim 
that Euclid could have developed proportionality in Book~I. 

\subsection{Bernay's Supplement II}
The eighth edition (1956) of Hilbert's {\em Foundations of Geometry} \cite{hilbert1899}  contained
Supplement II by Paul Bernays (which is still in the cited tenth edition).  The Supplement is 
entitled, {\em A simplified development of the theory of proportion}.   Bernay's definition 
of proportion, though phrased in terms of congruent angles instead of parallel lines,
is easily seen to be equivalent to our Definition~\ref{definition:proportion}.   Bernay sketches
a proof of the interchange theorem.

Bernays next considers what he calls the ``fundamental theorem of proportion'', namely,
if two parallels delineate the segments $a,a^\prime$ on one side of an angle and $b,b^\prime$
on the other side of the angle, then $a:a^\prime = b:b^\prime$.  
 Bernays proves the fundamental theorem (on pp.~204-205).
 The proof is both clear
and easily accessible, so we do not repeat it.  It depends on the 
theorem that the three bisectors of the angles of a triangle intersect in one point
(the incenter). Since we wish
to claim that Bernays's proof of the fundamental theorem uses only methods
from the part of Book I not using equal figures,  we need to check that 
the incenter theorem can be proved by those means.   That theorem
 occurs in Euclid, implicitly rather than explicitly, as Prop.~IV.4,
which is stated, {\em In a given triangle to inscribe a circle}.   But the first
half of Euclid's proof of IV.4 does not mention circles and proves the 
incenter theorem using only congruent triangles and perpendiculars. 
 
 Bernays also proves the existence and uniqueness of the fourth proportional. 
Hence Bernays's 1956 Supplement~II provides almost what we need in this paper:  a development
of the theory of proportions based on the methods of Euclid Book~I.   It falls short of that 
requirement only by needing a couple of simple theorems from Book~III about circles, which 
do not use ``equal figures'' for their proofs.  Nevertheless, we shall show
 (in \S\ref{section:kupffer2} below) how to 
overcome even this small defect, using a proof due to Kupffer. 

In summary:  Bernays's Supplement II proves all the theorems
listed in \S\ref{section:proportionality} as ``Theorems'', using
techniques to be discussed in more detail below, but generally 
acceptable for our purposes.  In the next section, we will 
prove the other results (lemmas and corollaries) from \S\ref{section:proportionality}.
 
\subsection{Kupffer's development of proportionality}\label{section:kupffer2}
As it turns out, Bernays was by no means the first ones to consider the possibility 
of developing the theory of proportions without using the theory of 
magnitudes (and Archimedes's axiom) as in Euclid Book~IV.   This was 
done by Karl Kupffer possibly among others,  as early as 1893.
We cite \cite{schur1902},  but that 1902 letter just calls attention
to  Kupffer's 1893 lecture, whose audience included Schur.
  Kupffer gave {\em two} elementary 
proofs of  Theorem~\ref{theorem:proportion-interchange}, both of which Euclid
``could have given.''%
\footnote{Kupffur gave two proofs of the interchange theorem in 1893. 
Bernays gave a proof in 1956, identical to Kupffur's first proof,
 that Bernays attributed to Federigo Enriques's 1911 book 
\cite{enriques1911}.  The relevant material is in a chapter written by someone
else, namely Giovanni Vailati.  On p.~239,  Vailati gives Kupffur's {\em second} proof,
with credit and citation,  and mentions his first proof,  but then says that the first proof is actually due 
to Weierstrass, who (Vailati says) was the first to develop proportion theory 
without the axiom of Archimedes.  But Vailati gives no citation to support this 
claim, and I could not pick up the trail.}

Since Schur's paper is in German, and since we want a detailed proof
to serve as a basis for formalization, we give both proofs here.
The first proof uses some theorems  which, as Euclid stated them,
mention circles, but are easily stated and proved without mentioning
circles, by elementary means from Book~I. Namely,
Prop.~III.21 (chords subtending the same arc are equal),
Prop.~IV.5 (three points determine a circle), and 
the following lemma, which is itself proved from III.21.
Prop.~III.21 uses III.20, which uses I.5 and I.32; the point is that
the use of equal figures starts with I.34,  so III.21 could be 
reached in two propositions after I.32, without using equal figures.
The proof of Prop.~IV.5 references only I.10 and I.4.  So these theorems
can all be proved without much of a detour from Book~I,  as we 
checked carefully.%
\footnote{\,Euclid~III.20 and III.21 need repairs both in 
statement and proof.  Some of these problems have been known
for centuries (see Heath's commentary \cite{euclid1956}). But
Heath does not remark on the step in Euclid's proof of III.20 for which 
Euclid gives no 
justification, but which is difficult to prove formally, and 
the final ``Therefore etc.'' obscures the difficult proof that the 
two cases Euclid presents (and the one he does not present)
are actually exhaustive. 
}    

\begin{lemma}[Cyclic quadrilateral theorem]\label{lemma:cyclicquadrilateral}
Let $ABCD$ be a convex quadrilateral whose diagonals meet at $O$.
Suppose angles $OAB$ and $ODC$ are equal.  Then the four vertices
lie on a circle.
\end{lemma}

 \begin{figure}[ht]
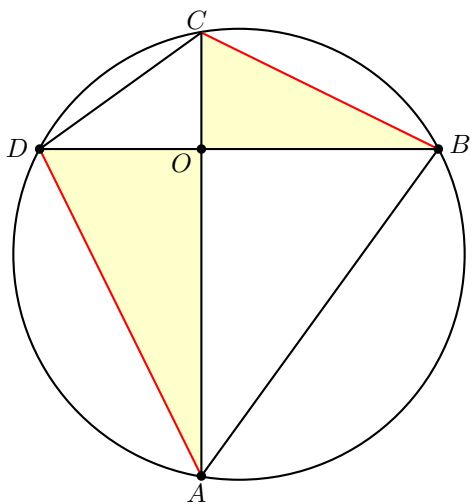

\caption{The cyclic quadrilateral theorem}
\label{figure:cyclicquadrilateral}
\FigureCyclicQuadrilateral
\end{figure}

\noindent{\em Proof}.  By Prop.~IV.5, 
any three non-collinear points lie on a circle, so let $K$ be a
circle containing $A$, $C$, and $D$.  Then point $O$ is inside $K$,
since it lies between $A$ and $C$.  Hence, by the line-circle axiom,
line $DO$ meets circle $K$ in two points; one of these points is $D$.
Call the other one $E$. Then angle $ODC$ and angle $OAE$ subtend
the same chord $EC$.  Hence by III.21, they are equal.  
But angle $ODC$ is equal to angle $OAB$ by hypothesis. Hence
angles $OAB$ and $OAE$ are equal.  Hence $B$ lies on line $AE$.
But $B$ also lies on line $OD$.  Hence $B$ is the intersection 
point of $DO$ and $AE$.  But that intersection point is $E$.
Hence $B=E$. Hence $B$ lies on circle $K$.  That completes the proof.

\subsection{Kupffer's first proof of the interchange theorem}

 Recall that the interchange theorem (Theorem~\ref{theorem:proportion-interchange}) 
 is:  if $r:s = p:q$, then $r:p = s:q$.
 
 \begin{figure}[ht]
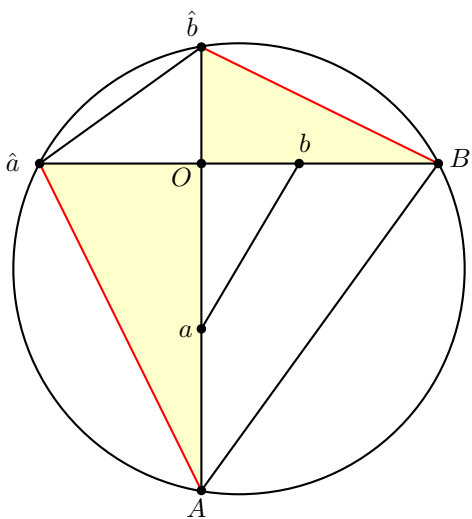

\caption{Kupffer's first proof of the interchange theorem}
\label{figure:kupffer1}
\FigureKupfferOne
\end{figure}

{\em Proof}.
Suppose $r:s = p:q$.  Then by Definition~\ref{definition:proportion},
there is a right angle $AOB$ with $AO = r$
and $BO=s$, and point $a$ on ray $OA$ and $b$ on ray $OB$ such that 
$Oa = p$ and $Ob = q$, and $AB\ ||\ ab$.  Without loss of generality 
we may assume $a < A$.  Let $\hat a$ and $\hat b$ be points on the other 
side of $O$ from $b$ and $a$, respectively, such that $O\hat a = Oa$
and $O\hat b = Ob$.  Angle $aOb$ is equal to angle $\hat a O \hat b$,
since they are vertical angles. 
 Then triangle $aOb$ is congruent to triangle $\hat ao\hat b$,
by SAS.  Since $AB\ ||\ ab$, angle $OAB$ is equal to angle $Oab$,
and hence also to angle $O\hat a \hat b$.  Then $A,B,\hat a, \hat b$
form a cyclic quadrilateral, i.e., all four lie on a circle, 
by Lemma~\ref{lemma:cyclicquadrilateral}. Then $OA\hat a$ and $OB\hat b$ 
have corresponding angles equal, since their angles at $O$ are vertical angles,
and their angles at $A$ and $B$ both subtend the same arc $\hat a \hat b$,
so they are equal by Euclid~III.21, and their angles at $\hat a$ and 
$\hat b$ subtend the same arc $AB$, so they are also equal by 
Euclid~III.21.   

  By Corollary~\ref{lemma:proportion1}, $OA:O\hat a = OB:O\hat b$.
But $O \hat a = Oa$ and $O \hat b = Ob$.  Therefore $AO:Oa = OB:Ob$;
that is $r:p = s:q$.  That completes Kupfer's proof.\
\FloatBarrier

\subsection{Kupffer's second proof of the interchange theorem}

This does not use any theorems about circles.  It depends only 
on the fact that the three altitudes of a triangle meet in a point (the ``orthocenter'').
That theorem was known to Archimedes, but it does not occur in Euclid.  
It does, however, have a short proof using the methods of Book~I, without
using equal figures, and hence 
certainly could have been proved by Euclid.  Such a proof 
can be found, for example, in \cite{hartshorne}, p.~54.

Pascal's theorem concerns three rays in a plane, meeting at a common point
$O$.  Kupffer's insight was that the special case when the three rays 
form two right angles is enough for the theory of proportionality, and 
that special case can be proved using existence of the orthocenter.

\begin{theorem}[Pascal, Kupffer's version] \label{theorem:pascal-schur}
Let $ABC$ be three distinct points on line $OA$ perpendicular at $O$ to 
the line containing three distinct points $A^\prime B^\prime C^\prime$ , all on the same side of $O$ and
distinct from $O$. 
Suppose that $AB^\prime\ ||\ BA^\prime$ and $BC^\prime\ ||\ CB^\prime$.  Then $AC^\prime\ ||\ CA^\prime$.
\end{theorem}

\begin{figure}
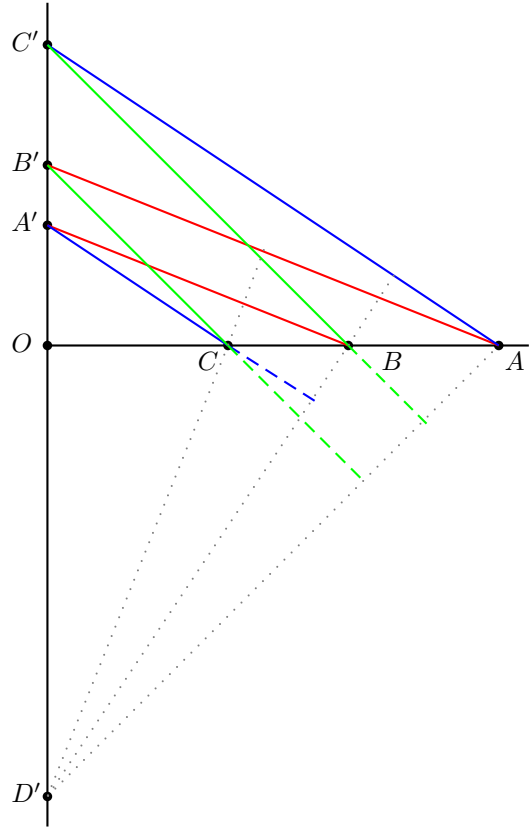

\caption{Kupffer's second proof.  Given two pairs of parallel lines (red and green),
prove that the blue lines are parallel too.}
\label{figure:kupffer2}
\FigureKupfferTwo
\end{figure}

\noindent{\em Remark}.  Pascal's theorem differs from this 
theorem in not requiring angle $AOC^\prime$ to be a right angle.
This theorem implies Pascal's theorem easily with the aid 
of Desargues's theorem.  On the other hand, it is known that 
Desargues's theorem can be proved with three applications of 
Pascal's theorem (but not this special case).  At any rate,
this version can be proved without Desargues.
\medskip

\noindent{\em Proof}. We assume known that the three altitudes of a triangle meet in a point,  the {\em orthocenter}
of the triangle.   We just 
translate Schur's proof from German in \cite{schur1902}, filling in no additional steps.
Construct the perpendicular line  from $B$ to $CA^\prime$.   Let it meet line $OD$ in point $D^\prime$.
Then $C$ is the orthocenter of the triangle $BA^\prime D$.  Thus $CD^\prime \perp BA^\prime$ and therefore
$CD^\prime \perp AB^\prime$.   Therefore $C$ is the orthocenter of triangle $AB^\prime D^\prime$.  Therefore
$AD^\prime \perp CB^\prime$ and also $AD^\prime \perp BC^\prime$.  Finally, $B$ is the orthocenter
of the triangle $AC^\prime D^\prime$.   Hence $AC^\prime \perp BD^\prime$.  Therefore $AC^\prime || CA^\prime$,
which is what was to be proved. 

\begin{corollary} If  $a:b = p:q$, then $a:p = b:q$.
\end{corollary}

{\em Remark}. This is a proof of Theorem~\ref{theorem:proportion-interchange} by 
the methods of Euclid Book~I, without even a slight detour into Book~III.
\medskip

\noindent{\em Proof}. Given $a:b= p:q$,
by Definition~\ref{definition:proportion},
there is a right angle with points $B^\prime$,
$A$, $A^\prime$,$B$ as in 
Fig.~\ref{figure:kupffer2}, with the red lines $B^\prime A$
and $A^\prime B$ parallel, and $a = OA$, $b = OB^\prime$, 
$p = OB$, and $q = OA^\prime$.  Now construct $C$ 
on ray $OB$ so that $OC = OB^\prime=b$, and 
construct $C^\prime$ on ray $OA^\prime$ such that 
$OC^\prime = OB = p$.  Then by
Theorem~\ref{theorem:pascal-schur}, the blue lines 
$C^\prime A$ and $A^\prime C$ are parallel. Then 
by Definition~~\ref{definition:proportion},
$OA:OC^\prime = OC:OA^\prime$.  That is, $a:p =b:q$.
That completes the proof.

\section{Conclusions}
We have given a definition of ``equal figures''  in the spirit of Euclid, using 
a diagram similar to the diagram for Prop.~I.44.   The fundamental 
properties of this defined notion seem to require (parts of) the theory
of proportion.   Our work then fell into two parts:

\begin{itemize}
\item  Using methods like those in Euclid Book~I,  as well as 
the elementary theory of proportions,  we proved  all the theorems 
of Book~I,  and all the equal-figures axioms used in  \cite{beeson2019},  using the defined notion of ``equal figures.''
\smallskip

\item We then showed, using theorems of Kupffer and Bernays, that 
the required theory of proportions can also be developed using the 
methods of Euclid Book~I.   In particular Desargues's theorem is not needed.
\smallskip

\item The final result is that it is possible to use the new definition
of ``equal figures'' to justify Euclid's use of that notion without
appealing to the Common Notions,  or to the ``equal figures'' 
axioms  used in \cite{hartshorne} or \cite{beeson2019}.
\smallskip

\item And Euclid could have done it. 
\end{itemize}

\section*{Appendix: Listing of the equal-figures axioms}
This formal listing is intended for reference.  Only ASCII
symbols are used, to facilitate cut-and-paste to computer-readable files.
Polish notation can be found in \cite{beeson2019}. 
\medskip

\noindent
{\em Congruent triangles are equal.} \vskip-5pt
\begin{lstlisting}[breaklines,breakatwhitespace,basicstyle=\ttfamily\footnotesize]
axiom-congruentequal 
forall A B C a b c, TC(A,B,C,a,b,c) ==> ET(A,B,C,a,b,c)
\end{lstlisting}

\noindent
{\em Triangles with the same vertices are equal.} \vskip-5pt
\begin{lstlisting}[breaklines,breakatwhitespace,basicstyle=\ttfamily\footnotesize]
axiom-ETpermutation
forall A B C a b c, ET(A,B,C,a,b,c) ==> ET(A,B,C,b,c,a) /\ ET(A,B,C,a,c,b) /\ ET(A,B,C,b,a,c) /\ ET(A,B,C,c,b,a) /\ ET(A,B,C,c,a,b)\end{lstlisting}

\noindent
{\em Triangle equality is a symmetric relation}\vskip-5pt
\begin{lstlisting}[breaklines,breakatwhitespace,basicstyle=\ttfamily\footnotesize]
axiom-ETsymmetric
forall A B C a b c, ET(A,B,C,a,b,c) ==> ET(a,b,c,A,B,C)\end{lstlisting}

\noindent
{\em Quadrilaterals with the same vertices are equal}\vskip-5pt
\begin{lstlisting}[breaklines,breakatwhitespace,basicstyle=\ttfamily\footnotesize]
axiom-EFpermutation
forall A B C D a b c d, EF(A,B,C,D,a,b,c,d) ==> EF(A,B,C,D,b,c,d,a) /\ EF(A,B,C,D,d,c,b,a) /\ EF(A,B,C,D,c,d,a,b) /\ EF(A,B,C,D,b,a,d,c) /\ EF(A,B,C,D,d,a,b,c) /\ EF(A,B,C,D,c,b,a,d) /\ EF(A,B,C,D,a,d,c,b)\end{lstlisting}
\goodbreak

\noindent
{\em Halves of equals are equal} \vskip-5pt
\begin{lstlisting}[breaklines,breakatwhitespace,basicstyle=\ttfamily\footnotesize]
axiom-halvesofequals
forall A B C D a b c d, ET(A,B,C,B,C,D) /\ OS(A,B,C,D) /\ ET(a,b,c,b,c,d) /\ OS(a,b,c,d) /\ EF(A,B,D,C,a,b,d,c) ==> ET(A,B,C,a,b,c)\end{lstlisting}

\noindent{\em Equal quadrilaterals is a symmetric relation} \vskip-5pt
\begin{lstlisting}[breaklines,breakatwhitespace,basicstyle=\ttfamily\footnotesize]
axiom-EFsymmetric
forall A B C D a b c d, EF(A,B,C,D,a,b,c,d) ==> EF(a,b,c,d,A,B,C,D)\end{lstlisting}
\goodbreak

\noindent{\em Equal quadrilaterals is a transitive relation} \vskip-5pt
\begin{lstlisting}[breaklines,breakatwhitespace,basicstyle=\ttfamily\footnotesize]
axiom-EFtransitive
forall A B C D P Q R S a b c d, EF(A,B,C,D,a,b,c,d) /\ EF(a,b,c,d,P,Q,R,S) ==> EF(A,B,C,D,P,Q,R,S)\end{lstlisting}
\goodbreak

\noindent{\em Equal triangles is a transitive relation} \vskip-5pt
\begin{lstlisting}[breaklines,breakatwhitespace,basicstyle=\ttfamily\footnotesize]
axiom-ETtransitive
forall A B C P Q R a b c, ET(A,B,C,a,b,c) /\ ET(a,b,c,P,Q,R) ==> ET(A,B,C,P,Q,R)\end{lstlisting}
\goodbreak

\noindent{\em Cutting off equal triangles from equal triangles
yields equal quadrilaterals} \vskip-5pt
\begin{lstlisting}[breaklines,breakatwhitespace,basicstyle=\ttfamily\footnotesize]
axiom-cutoff1
forall A B C D E a b c d e, BE(A,B,C) /\ BE(a,b,c) /\ BE(E,D,C) /\ BE(e,d,c) /\ ET(B,C,D,b,c,d) /\ ET(A,C,E,a,c,e) ==> EF(A,B,D,E,a,b,d,e)\end{lstlisting}
\goodbreak

\noindent{\em Cutting off equal triangles from equal quadrilaterals
yields equal quadrilaterals} \vskip-5pt
\begin{lstlisting}[breaklines,breakatwhitespace,basicstyle=\ttfamily\footnotesize]
axiom-cutoff2
forall A B C D E a b c d e, BE(B,C,D) /\ BE(b,c,d) /\ ET(C,D,E,c,d,e) /\ EF(A,B,D,E,a,b,d,e) ==> EF(A,B,C,E,a,b,c,e)\end{lstlisting}
\goodbreak

\noindent{\em Pasting equal triangles yields equal triangles} \vskip-5pt
\begin{lstlisting}[breaklines,breakatwhitespace,basicstyle=\ttfamily\footnotesize]
axiom-paste1
forall A B C D E a b c d e, BE(A,B,C) /\ BE(a,b,c) /\ BE(E,D,C) /\ BE(e,d,c) /\ ET(B,C,D,b,c,d) /\ EF(A,B,D,E,a,b,d,e) ==> ET(A,C,E,a,c,e)\end{lstlisting}
\goodbreak

\noindent{\em Cutting off a triangle makes an unequal triangle} \vskip-5pt
\begin{lstlisting}[breaklines,breakatwhitespace,basicstyle=\ttfamily\footnotesize]
axiom-deZolt1
forall B C D E, BE(B,E,D) ==> ~ET(D,B,C,E,B,C)\end{lstlisting}
\goodbreak

\noindent{\em Cutting off a quadrilateral makes an unequal triangle} \vskip-5pt
\begin{lstlisting}[breaklines,breakatwhitespace,basicstyle=\ttfamily\footnotesize]
axiom-deZolt2
forall A B C E F, TR(A,B,C) /\ BE(B,E,A) /\ BE(B,F,C) ==> ~ET(A,B,C,E,B,F)\end{lstlisting}
\goodbreak

\noindent{\em Pasting equal triangles to equal quadrilaterals yields equal quadrilaterals} \vskip-5pt
\begin{lstlisting}[breaklines,breakatwhitespace,basicstyle=\ttfamily\footnotesize]
axiom-paste2
forall A B C D E M a b c d e m, BE(B,C,D) /\ BE(b,c,d) /\ ET(C,D,E,c,d,e) /\ EF(A,B,C,E,a,b,c,e) /\ BE(A,M,D) /\ BE(B,M,E) /\ BE(a,m,d) /\ BE(b,m,e) ==> EF(A,B,D,E,a,b,d,e)\end{lstlisting}
\goodbreak

\noindent{\em Pasting equal triangles to equal triangles yields equal quadrilaterals} \vskip-5pt
\begin{lstlisting}[breaklines,breakatwhitespace,basicstyle=\ttfamily\footnotesize]
axiom-paste3
forall A B C D M a b c d m, ET(A,B,C,a,b,c) /\ ET(A,B,D,a,b,d) /\ BE(C,M,D) /\ BE(A,M,B) \/ EQ(A,M) \/ EQ(M,B) /\ BE(c,m,d) /\ BE(a,m,b) \/ EQ(a,m) \/ EQ(m,b) ==> EF(A,C,B,D,a,c,b,d)\end{lstlisting}
\goodbreak

\noindent{\em Pasting equal quadrilaterals yields equal quadrilaterals} \vskip-5pt
\begin{lstlisting}[breaklines,breakatwhitespace,basicstyle=\ttfamily\footnotesize]
axiom-paste4
forall A B C D F G H J K L M P e m, EF(A,B,m,D,F,K,H,G) /\ EF(D,B,e,C,G,H,M,L) /\ BE(A,P,C) /\ BE(B,P,D) /\ BE(K,H,M) /\ BE(F,G,L) /\ BE(B,m,D) /\ BE(B,e,C) /\ BE(F,J,M) /\ BE(K,J,L) ==> EF(A,B,C,D,F,K,M,L)\end{lstlisting}
 
\section*{Appendix: Where the equal-figures axioms are used}
The following listing shows all the lines in the formal 
development of \cite{beeson2019} that are justified by the 
equal-figure axioms other than the axioms {\tt ETpermutation},
{\tt EFpermutation}, and the axioms asserting that {\tt ET} and 
{\tt EF} are equivalence relations.  The middle entry in each 
line is the statement justified; in most cases, the reader will 
be able to identify the corresponding line in Euclid's own proof, which
will either be justified by a common notion, or not justified at all. 
To decode the statements: for example in {\tt EFADGBFEGC}, 
the initial {\tt EF} means ``equal figures'' and the statement means
that $ADGB$ and $FEGC$ are equal quadrilaterals.

\medskip

\begin{verbatim}
Prop35A.prf:      EFADGBFEGC     axiom:cutoff1
Prop43.prf:       EFAKGBAKFD     axiom:cutoff1
Prop43.prf:       EFGBEKFDHK     axiom:cutoff2
EFreflexive.prf:  EFabcdabcd     axiom:paste3
Prop35A.prf:      EFADCBFEBC     axiom:paste2
Prop35A.prf:      EFCDABBEFC     axiom:paste2
Prop35A.prf:      EFBAECCFDB     axiom:paste3
Prop42.prf:       EFABECFECG     axiom:paste3
Prop42B.prf:      EFABECabec     axiom:paste3
Prop45.prf:       EFABCDFKML     axiom:paste4
Prop47B.prf:      EFFBAGDBML     axiom:paste3
Prop48.prf:       EFBCEDBced     lemma:paste5
paste5.prf:       EFDBCLdbcl     axiom:paste2
paste5.prf:       EFBDECbdec     axiom:paste2
squaresequal.prf: EFBADCbadc     axiom:paste3
Prop48.prf:       EFACKHAckh     lemma:squaresequal
Prop48.prf:       EFABFGABfg     lemma:squaresequal
Prop39A.prf:      NOETDBCEBC     axiom:deZolt1
Prop39A.prf:      NOETEBCDBC     axiom:deZolt1
Prop48A.prf:      NOETDABFAE     axiom:deZolt2
Prop48A.prf:      NOETdabfae     axiom:deZolt2
Prop37.prf:       ETCBACBD       axiom:halvesofequals
Prop38.prf:       ETEFDCBA       axiom:halvesofequals
Prop48A.prf:      ETABDabd       axiom:halvesofequals
paste5.prf:       ETMCLmcl       axiom:halvesofequals
paste5.prf:       ETECLecl       axiom:halvesofequals
\end{verbatim}
 
\bibliographystyle{plain-annote}

\end{document}